\documentclass[12pt]{article}
\usepackage{times,amssymb,amsmath}

\begin{document}
\newtheorem{theorem}{Theorem}[section]
\newtheorem{proposition}[theorem]{Proposition}
\newtheorem{lemma}[theorem]{Lemma}
\newtheorem{corollary}[theorem]{Corollary}
\newtheorem{definition}{Definition}[section]

\newcommand{\1}{{\bf 1}}
\newcommand{\A}{{\mathcal A}}
\newcommand{\D}{{\mathbb D}}
\newcommand{\E}{{\mathbb E}}
\newcommand{\F}{{\mathcal F}}
\newcommand{\PP}{{\mathbb P}}
\newcommand{\R}{{\mathbb R}}
\newcommand{\RR}{{\mathcal R}}
\newcommand{\pathR}{{\mathcal{\rm I\!R}}}
\newcommand{\I}{{\mathcal I}}
\newcommand{\LL}{{\mathcal L}}
\newcommand{\na}{{\bf \nabla}}
\newcommand{\kk}{{\mathfrak k}}
\newcommand{\HH}{{\mathcal H}}
\def\g{\mathfrak g}
\def\n#1{|\kern-0.24em|\kern-0.24em|#1|\kern-0.24em|\kern-0.24em|}
\def\paral{/\kern-0.55ex/}
\def\parals_#1{/\kern-0.55ex/_{\!#1}}
\renewcommand{\thefootnote}{}

\title{Some families of $q$-vector fields  on path spaces}
\author{K. D. Elworthy and Xue-Mei Li}
\date{}
\maketitle

\noindent
\footnote{ key words: path space, differential forms, Riemannian manifolds,
heat equation, integration by parts, Bismut formula, Malliavin calculus.}
\noindent\footnote{Research partially supported by NSF  research grant
DMS 0072387, EPSRC GR/NOO 845, the Alexander von Humboldt stiftung.}
\begin{abstract}
Some families of $H$-valued vector fields with calculable Lie brackets are given.
These provide examples of vector fields on path spaces with a divergence
and we show that versions of Bismut type formulae for forms on a compact
Riemannian manifold arise as projections of the infinite dimensional theory.
\end{abstract}

Let $M$ be a compact Riemannian manifold and $\{P_t:t\ge 0\}$ its
heat semigroup acting on differential forms. Thus if $\phi$ is a
bounded measurable (or square integrable) form on $M$ we have
$$\left\{
 \begin{array}{lll}
{\partial \over \partial t}(P_t\phi)&=&{1\over 2}\Delta(P_t\phi),
\hskip 18pt t>0\\
P_0\phi&=&\phi\end{array}\right.$$
where we use the sign convention
$$\Delta=-(d\delta+\delta d).$$
The formula for the exterior derivative of $P_t\phi$
\begin{equation}
\label{1}
d(P_t\phi)(V_0)=
{1\over t}\E \phi\left(
W_t^{(q)}\int_0^t (W_s^{(q)})^{-1} \iota_{\langle -,dx_s\rangle}
 W_s^{(q+1)}(V_0)\right),
\end{equation}
 where $\phi$ is a bounded measurable $q$-form and
 $V_0\in\wedge^{q+1}T_{x_0}M$ was given in \cite{Elworthy-Li-forms-CR},
 extending Bismut's formula for the special case of $q=0$, \cite{Bismut-book}
\begin{equation}
\label{2}
d(P_tf)(V_0)={1\over t}
\E f(x_t)\int_0^t\langle W_s^{(1)}(V_0), dx_s\rangle
\end{equation}
for $f:M\to \R$ bounded and measurable and $V_0\in T_{x_0}M$. In
these formulae the expectation $\E$ is with respect to a  Brownian
motion $\{x_t: 0\le t<\infty\}$ on $M$ starting from $x_0$ and for
$q=1,2,\dots$ we use the `damped' parallel translations
 $W_s^{(q)}:\wedge^qT_{x_0}M\to \wedge^q T_{x_s}M$, $0\le s<\infty$ of
 $q$-vectors  along the sample paths of $(x_\cdot)$ given by the covariant
 equation along $(x_\cdot)$:
\begin{equation}\left\{
\begin{split}
\label{3}
{D\over \partial s}(W_s^{(q)}(V_0))
&=-{1\over 2}{\RR}_{x_0}^q(W_s^{(q)}(V_0))\\
W_0^{(q)}(V_0)&=V_0,
\end{split}\right.
\end{equation}
where $\RR^q_x:\wedge^qT_xM\to \wedge^q T_xM$ is the q-th Weitzenbock
curvature defined by
$$\Delta \phi=\nabla^*\nabla \phi-\phi(\RR^q_\cdot -)$$
for $\phi$ a smooth $q$-form. See, for example, \cite{Airault76},
\cite{Elworthy-book}, \cite{Ikeda-Watanabe} and \cite{Elworthy-Stflour}.
In particular $\RR^1$ is the Ricci curvature $Ric^\#:TM\to TM$, i.e.
$ \langle \RR_x^1 v,u\rangle =Ric(v,u)$,  $u, v\in T_xM$.

Formula (\ref{1}) has been refined and extended by Driver and Thalmaier
\cite{Driver-Thalmaier-2001}, giving analogous results
on various operators on vector bundles over$M$, e.g. the square of
the Dirac operator. Rather different types of Bismut type formulae
for covariant derivatives of operators on vector bundles were
obtained by Norris in  \cite{Norris93} by different methods. The
original proof of (\ref{1}) in \cite{Elworthy-Li-forms-CR} was
to derive it by applying the method of conditional expectation
from \cite{Elworthy-Yor} to the earlier, non-intrinsic, formula of
Li, \cite{thesis}, \cite{Elworthy-Li-JFA}, see also \cite{Elworthy-flow},
\begin{equation}
\label{4}
d(P_t\phi)_{x_0}={1\over t}\E\int_0^t \langle T\xi_s-, dx_s\rangle \wedge
\xi_t^*(\phi)
\end{equation}
for $\phi$ a bounded measurable $q$-form. Here $\{\xi_s: 0\le s<\infty\}$ is
a gradient Brownian flow on $M$ (see below), $x_t=\xi_t(x_0)$ for $x_0\in M$,
$T\xi_t:TM\to TM$ is the (random) derivative of $\xi_t$,
and $\xi_t^*\phi$ is the pull back of $\phi$:
\begin{eqnarray*}
\xi_t^*\phi(V):&=&\phi\left(\wedge^q(T\xi_t)(V)\right)\\
&=&\phi\left(T\xi_t(v^1)\wedge\dots\wedge T\xi_t(v^q)\right),
\hskip 15pt V=v^1\wedge\dots\wedge v^q\in \wedge^qTM
\end{eqnarray*}
In fact the same proof allows the use of any Brownian flow which
has Levi-Civita connection as its LeJan-Watanabe connection in the
sense of \cite{Elworthy-LeJan-Li-book}. All the formulae can be extended to
obtain differentiation formulae for the corresponding heat kernels
and this is done in the references cited.

  Bismut's formula (\ref{2}) can be obtained from infinite dimensional
integration by parts formulae by considering the cylindrical
functions $F(\sigma_\cdot)=f(\sigma_t)$ on the  space of paths
over $M$ and indeed Driver's integration by parts formula can be
derived from it, as described in \cite{Elworthy-Li-Ito} following
comments by Nualart. Since integration by parts theorems for forms
on path spaces are not yet well understood it is interesting to
ask if analogous results hold for (\ref{1}) and (\ref{4}) in the
context of the $L^2$ theory of differential forms being developed
in \cite{Elworthy-Li-Hodge-1} and \cite{Elworthy-Li-Hodge-2} and
to derive (\ref{1}) and ({\ref{4}) by the methods used there. Our
approach is very much in the spirit of Bismut's approach to
Malliavin Calculus \cite{Bismut-Durham} and of Fang-Franchi
for Lie groups \cite{Fang-Franchi97}. We do not attack the two
challenges of :
(i) deriving integration by parts formulae for forms on path spaces from
(\ref{1}) and ({\ref{4}), as done for functions in \cite{Elworthy-Li-Ito};
and (ii) deriving the more general results of Driver and Thalmaier
 by similar methods.

 In \S1 below we define the class of flows used in (\ref{1}). In \S2 the
 infinite dimensional theory in \cite{Elworthy-Li-Hodge-1},
 \cite{Elworthy-Li-Hodge-2} is briefly  described, and in \S3 the two parts
 are combined to both give a proof of (\ref{1}) and ({\ref{2}) and to give
 some interesting examples of q-vector fields on path spaces which have
 explicitly calculable divergences. On the way we obtain, from our stochastic
 differential equation (\ref{5}), a family of H-vector fields on Wiener space,
 which for gradient stochastic differential equations form a commuting family
 and in certain other cases have easily computable Lie brackets.
 This gives a $q+1$-vector field $\tilde V_0$
on the path space $C_{x_0}M$ of $M$ which has a divergence with respect
to the Wiener measure $\mu_{x_0}$ in the sense that the divergence
 $div\tilde V_0$ is a  vector field on $C_{x_0}M$ satisfying equation
(\ref{14}) below and  such that (\ref{1}) can be
written:
$$\int_{C_{x_0}M} d\phi^t(\tilde V_0)\; d\mu_{x_0}
  =  - \int_{C_{x_0}M}\phi^t(div\tilde V_0)\;d\mu_{x_0}$$
for $\phi^t$ the cylindrical $q$-form on $C_{x_0}M$ obtained from
a differential form $\phi$  on the manifold $M$:
$$\phi^t(V)=\phi(x_t)(V_{t,\dots,t}), \hskip 15pt
\hbox{for $V$ any $q$-vector field}.$$

Various versions of formula (\ref{2}) are derived from this. For invariant
stochastic differential equations on Lie groups no conditioning is needed
in (\ref{4}), c.f. \cite{Fang-Franchi97}, and we obtain
(\ref{formula-Li-group}). Using gradient systems we give
in (\ref{Bismut-formula}) the minor extension to (\ref{4}) mentioned in
\cite{Elworthy-Li-forms-CR} to the case of semigroups on forms with
generators of the form ${1\over 2}\Delta+\LL_A$. We also consider
the connections used by Ikeda and Watanabe in \cite{Ikeda-Watanabe},
sometimes known as Riemann-Cartan-Weyl connections
 \cite{Rapoport-Sternberg84},  where
the formulae (\ref{2}) and (\ref{4}) have a pleasant form
(\ref{41}), although the infinitesimal generators of the semigroup on
forms are rather complicated. However the main aims of the article
are to study some interesting classes of $q$-vector fields on $C_{x_0}M$
and to show that some finite dimensional formulae can be considered as
projections of the infinite dimensional theory.

\section{Stochastic flows and LeJan-Watanabe connections}
\label{section-1}
Consider a Stratonovich stochastic differential equation
\begin{equation}
\label{5}
dx_t=X(x_t)\circ dB_t +A(x_t)dt
\end{equation}
on $M$ driven by a $m$-dimensional Brownian motion $\{B_t, 0\le t<\infty\}$.
Here $X: M\times \R^m\to TM$ and $A: M\to TM$ are assumed to be smooth and
$X(x):= X(x,-):\R^m\to T_xM$ is assumed to be linear, surjective for each $x$
and to induce the inner product $\langle, \rangle_x$ on $T_xM$ given
by the Riemannian structure of $M$. This implies the Markov process of
solutions to (\ref{5}) has infinitesimal generator of the form
\begin{equation}
\label{6}
\A={1\over 2}\Delta +\LL_{Z}
\end{equation}
where $\LL_Z$ is Lie differentiation in the direction of some vector
field $Z$:
$$\LL_Z(f)(x)=df(Z(x)), \hskip 20pt x\in M.$$
The solution to (\ref{5}) from $x$ shall be denoted by
$\{\xi_t(x), 0\le t<\infty\}$. Write $x_t=\xi_t(x_0)$. We can, and will,
take versions which makes $\xi_t(x)$ continuous in $t$ and a $C^\infty$
diffeomorphism of $M$ in $x$.
\bigskip

Let $Y:TM\to \R^m$  be the adjoint of $X$, i.e. $Y_x=X(x)^*$, so that
$Y$ is the right inverse to $X$. As described in
 \cite{Elworthy-LeJan-Li-Tani}, or more generally in
 \cite{Elworthy-LeJan-Li-book}, $X$ induces a metric connection
 $\breve \nabla$ on $M$ defined by
\begin{equation}\label{7}
\breve \nabla_v U=X(x)d [x\mapsto Y_x(U(x))](v), \hskip 20pt v\in
T_xM\end{equation} for any smooth vector field $U$ on $M$. In
general this, which we call the {\it LeJan-Watanabe connection } of
our stochastic differential equation, has torsion. The torsion is given by
\begin{equation}
\label{8}\breve T(u,v)=X(x)dY(u,v),
\hskip 20 pt u.v\in T_xM
\end{equation}
where $dY$ refers to the exterior derivative of $Y$ considered as
an $\R^m$-valued differential 1-form on $M$, see Proposition 2.2.3
in \cite{Elworthy-LeJan-Li-book}. When the torsion vanishes $\breve \nabla$
 is  just the Levi-Civita connection, which we will always refer to
using $\nabla$, ${D\over \partial t}$ etc. There are two principal
classes of examples, Example 1.1 and Example 1.2 below, for which
 this holds.

\subsubsection*{Example 1.1. Gradient system.}
Here we consider an isometric immersion  $\alpha:M\to \R^m$ for
some $m$, e.g. by using Nash's theorem. Let $X(x): \R^m\to T_xM$
be the orthogonal projection identifying $T_xM$ with its image in
$\R^m$ under $d\alpha$. Then $Y_x: T_xM\to \R^m$ is just the
inclusion $(d\alpha)_x$ and (\ref{7}) is the classical formula for
the Levi-Civita connection of a sub-manifold of $\R^m$.

\subsubsection*{Example 1.2. Riemannian symmetric spaces.}
Let $M$ be a symmetric space $(K,H,\sigma)$. In particular $K$ is a Lie group
 acting transitively on $M$ and $H$ can be identified with the isotropy
 subgroup fixing the point $x_0$ of $M$, so that $k\mapsto kx_0$ gives a
diffeomorphism of $K/H$ with $M$. Assume that the Lie algebra
$\kk$ of $K$ has an inner product $\langle, \rangle_{\kk}$, say,
invariant under the adjoint action of $K$ on $\kk$. The action induces
a stochastic differential equation
$$X: M\times \kk \to TM$$ with
$X(x)e={d\over dt}[(\exp{(te)})\cdot x]_{t=0}$. It is surjective and so
induces a Riemannian structure $(\langle, \rangle_x, x\in M)$ on
$M$ which is $K$-invariant. The LeJan-Watanabe connection $\breve \nabla$
is the Levi-Civita connection for the Riemannian structure, as is proved in
 Corollary 1.4.9 of \cite{Elworthy-LeJan-Li-book}.
 \bigskip

An important class of examples which give connections other than the
Levi-Civita connections are the left invariant, and the right invariant
stochastic differential equations  on Lie groups. For $M$ a
Lie group with left invariant metric take the Lie algebra $\g$
as $\R^m$ and define $$X^L: G\times \g \to TG$$
to be left invariant with $X^L(e)\alpha=\alpha$, for $e$ the identity element
of $G$, i.e.
$$X^L(g)\alpha=TL_g(\alpha)$$
for $L_g: G\to G$ the left translation $x\mapsto g\cdot x$. Similarly if
 $R_g$ denotes right translation, $x\mapsto x\cdot g$, define
 $X^R:G\times \g \to TG$ by $X^R(g)\alpha=TR_g(\alpha)$.
 It is easy to see that the associated  connections in the sense of
 (\ref{7}) are just the canonical left and right  invariant connections,
  $\nabla^L$ and $\nabla^R$ say, of $G$, respectively. The flow $\xi_t^L$
 for the left invariant stochastic differential equation
 $dx_t=X^L(x_t)\circ dB_t$, where $B_t$ is a $\g$-valued Brownian motion,
is given by   $$\xi_t^L(x)=xg_t, \hskip 20 pt x \in G, t\ge 0.$$
 Here $g_t$ is the solution starting from the identity $e$ to
 $dg_t=X^L(g_t)\circ dB_t$. Similarly the solution $\xi_t^R(x)$,
  with $\xi_0^R(x)=x$  to the right invariant stochastic differential equation
$dx_t=X^R(x_t)\circ dB_t$ is given by  the right translation of $h_t$,
the solution starting from $e$. If the metric on $G$ is bi-invariant  both
stochastic differential equations, with $A\equiv 0$, have Brownian motions
as solutions. Also in that case the connections are torsion skew symmetric,
with $\nabla^L$ and $\nabla^R$ adjoint to each other. In particular the
torsions $T^L$ and $T^R$ are given by
\begin{equation}\label{torsion-of-Lie-groups}
T^L(u,v)\equiv -[U,V](g), \hskip 15pt  T^R(u,v)=[U,V](g),
\end{equation}  where $U$, $V$
are left invariant vector fields with $u=U(g)$ and $v=V(g)$.

In the bi-invariant case we can also treat $G$ as a symmetric space
by taking $K=G\times G$ with action $(g_1,g_2)x=g_1xg_2^{-1}$,
with $H=\Delta G\equiv\{ (g,g)\in K: g\in G\}$ the isotropy subgroup
fixing the identity $e$ of $G$. The symmetry map
 $\sigma:\g\times\g\to \g$ is given by
 $\sigma(\alpha, \beta)=(\beta, \alpha)$.
The stochastic differential equation induced by the action as described
in Example 1.2  above is just
\begin{equation}\label{sde-bi-invariant-Lie-group}
dx_t=TR_{x_t}\circ dB_t -TL_{x_t}\circ dB_t^\prime
\end{equation}
where $(B_t)$ and $(B_t^\prime)$ are two independent
Brownian motions on $\g$. The solution flow is then given by
$\xi_t(x)=g_tx(g_t^\prime)^{-1}$ where $g_t$ and $g_t^\prime$ are independent
 Brownian motions on $G$ and solutions to  the following right invariant
stochastic differential equations respectively:
$$dg_t=TR_{g_t}\circ dB_t, \hskip 18pt g_0=e,$$
$$dg_t^\prime=TR_{g_t^\prime}\circ dB_t^\prime, \hskip 18pt g_0=e,$$
The solutions $\xi_t(x)$, when $A\equiv 0$,  are also Brownian motions on $G$.

For more general homogeneous spaces $M=K/H$, if the
metric on $M$ is induced from an ad$_K$-invariant inner product
on $\kk$ and the orthogonal complement of the ker $X(x_0)$ is
ad$_H$ invariant then $\breve \nabla$, the connection associated to
the  stochastic differential equation defined as in Example 1.2 above,
is $K$-invariant. (see \cite{Elworthy-LeJan-Li-book}, Proposition 1.2.9).
Again in that case $$\xi_t(x)=k_t x, \hskip 18pt x\in M$$
for $k_t$ the solution to the equation $dk_t=TR_{k_t}\circ dB_t$ on $K$
 with initial point $k_0=e$.

\bigskip

\noindent{\bf Stochastic flows and differential forms.}
For $\{\xi_t: t\ge 0\}$ the flow of our stochastic differential equation
(\ref{5}) there are semigroups on differential forms defined by
$P_t\phi=\E\xi_t^*\phi$.
On $q$-forms the semigroup has infinitesimal generator $\A^{q}$
given on smooth forms by
$$\A^{q} \phi={1\over 2}\sum_{j=1}^m \LL_{X^j}\LL_{X^j}\phi+\LL_{A}(\phi)$$
for $\LL_{X^j}$ denoting the Lie differentiation in the direction $X^j$, etc
and with $X^j$ denoting the vector fields $X(\cdot)e^j$ for $e^1,\dots e^m$
the standard basis of $\R^m$. This comes from It\^o's formula, see
\cite{Elworthy-flow} or \cite{Elworthy-Stflour}. The generator can be written
\begin{equation}\label{11}
\A^q(\phi)=-(d\hat \delta+\hat \delta d)\phi+\LL_A(\phi)
\end{equation}
where, for $V\in \wedge^{q-1}TM$,
\begin{equation}
\label{hat-delta}
\begin{split}
\hat \delta \phi(V)=&-\sum_{j=1}^m(\iota_{X^j}\LL_{X^j}\phi)(V)\\
=&-trace \hat\nabla_-\phi(-,V),
\end{split}
\end{equation}
and $\hat \nabla$ is the adjoint connection to $\breve \nabla$,
in the sense of Driver \cite{Driver92} so the torsion of $\hat\nabla$
is $-\breve T$. See \cite{Elworthy-LeJan-Li-book}. This shows  $\A^q$
depends only on the connection $\breve \nabla$, the vector field $A$,  and  the Riemannian  metric on $TM$. In particular if $A\equiv 0$ in (\ref{5}) and  the associated
connection $\breve \nabla$ is the Levi-Civita connection we get  the usual
Hodge-Kodaira Laplacian, c.f. \cite{Kusuoka-degree}.

\section{Spaces of forms on the path space of $M$}

Fix $T>0$ and $x_0\in M$ and let $C_{x_0}M$ denote the Banach manifold of
continuous paths $\sigma:[0, T]\to M$ with $\sigma(0)=x_0$. Let $\mu_{x_0}$
be the measure induced on it by the  solution to (\ref{5}) starting
from $x_0$.

We can identify the tangent space $T_\sigma C_{x_0}M$ to $C_{x_0}M$ at
$\sigma$ with the space of continuous $v:[0,T]\to TM$ with $v(0)=0$ such that
$v(t)\in T_{\sigma(t)}M$ for each $0\le t \le T$. For $q=1,2,\dots$, let
$\wedge^qT_\sigma C_{x_0}M$ be the space of anti-symmetric $q$-vectors at
$x_0$ completed using the largest reasonable cross norm, so that its dual
space can be identified with the Banach space of continuous alternating
$q$-linear maps
$\phi: T_\sigma C_{x_0}M\times \dots \times T_\sigma C_{x_0}M\to \R$.
Thus $q$-forms on $C_{x_0}M$ are sections of the dual bundle
$(\wedge^qT_\sigma C_{x_0}M)^*$, e.g. see \cite{Lang-book}.

 As an example of a $q$-vector on $C_{x_0}M$ consider $\wedge^q(T_{x_0}\xi)(U)$
evaluated on a fixed $\omega  \in \Omega$, where $U\in \wedge^q T_{x_0}M$. It
is then a $q$-vector at $\{x_t(\omega): 0\le t\le T\}$. For $q=1$ it is the
tangent vector $\{T_{x_0}\xi_t(\cdot, \omega)(U): 0\le t \le T\}$ to
$C_{x_0}M$ at $x_\cdot(\omega)$. For higher $q$ we identify $q$-vectors $V$
with certain continuous maps $V:[0,T]^q \to \otimes^q TM$ such that
$V_{t_1,t_2,\dots, t_q}\in T_{\sigma(t_1)}\otimes \dots \otimes T_{\sigma(t_q)}M$.
Then, generalizing to allow $U$ to be in $\wedge^{q+1}C([0,T];T_{x_0}M)$,
$$\wedge^{q+1}(T_{x_0}\xi)(U)_{t_1,\dots, t_{q+1}}
=\left(T_{x_0}\xi_{t_1}\otimes \dots\otimes T_{x_0}\xi_{t_{q+1}}\right)
\left(U_{t_1,\dots,t_{q+1}}\right).$$
As $\omega$ varies in $\Omega$ such $q+1$ vectors do not strictly speaking
form $q+1$ vector fields on $C_{x_0}M$. However we can obtain $q$-vector
fields from them by conditioning, `filtering out the redundant noise',
to give
$\sigma\mapsto \E\{\wedge^{q+1}(T_{x_0}\xi)(U)\left|
 \xi_\cdot(x_0, \omega)=\sigma \right.\}$.
As shorthand, we write:
\begin{equation}
\overline {\wedge^{q+1}(T_{x_0}\xi)}_\sigma (U)
 =\E\{\wedge^{q+1}(T_{x_0}\xi)(U)\left|
 \xi_\cdot(x_0, \omega)=\sigma \right.\}.
 \end{equation}
  See \cite{Elworthy-Yor},\cite{Elworthy-LeJan-Li-book}. Below we shall often
  use  $\overline f$ denote the conditional expectation of a random function
  $f$ with respect to the $\sigma$-algebra generated by $\{x_t: 0\le t<\infty\}$.

An important class of forms on $C_{x_0}M$ are the cylindrical forms. In
particular if $\phi$ is a $q$-form on $M$ and $t\in [0, T]$ we have
the cylindrical $q$-form $\phi^{t}$ on $C_{x_0}M$ given by
$\phi^t(V)=\phi(V_{t,\dots, t})$.

For $\psi$ a $(q-1)$-form on $M$ we see, with $P_t$ as in the last section,
\begin{equation}
P_t(d\psi)=\E(d\psi)\left(\wedge^q(T_x\xi_t)U\right)
=\int_{C_{x_0}M} (d\psi)^t(\overline{\wedge^q T_{x}\xi})(U))\,
 d\mu_{x_0}(\sigma).
\end{equation}
In general we say that an $L^2$ $q$-vector field $V$ {\it has a divergence}
if there is an $L^2$ $(q-1)$-vector field $div V$ on $C_{x_0}M$ such that
for all smooth cylindrical $(q-1)$-forms $\phi$ on $C_{x_0}M$ we
have
\begin{equation}
\label{14}
\int_{C_{x_0}M} d\phi(V)\, d\mu_{x_0}
=- \int_{C_{x_0}M} \phi(div V) \, d\mu_{x_0}.
\end{equation}

For $U=u^1\wedge\dots \wedge u^q$, a primitive vector,  the non-intrinsic Bismut
 type formula (\ref{4}) will follow from showing that
$\overline{\wedge^q(T_{x_0}\xi)}(U)$ has a divergence with
\begin{eqnarray*}
&&div\left(\overline{\wedge^qT_{x}\xi}_\sigma(U)\right)_{t,\dots, t}\\
&&={1\over t}\E\left\{
\sum_{j=1}^q (-1)^{j}\left(\int_0^t \langle T\xi_s (u^j), dx_s\rangle
\right) \wedge^{q-1}(T\xi_t)(u^1\wedge\dots \widehat{u^j}\dots \wedge u^q)
\right.\\
&& \hskip 30pt \left. \left\vert \; \xi_\cdot(x_\cdot)=\sigma \right. \right\}
\end{eqnarray*}
at least for a gradient stochastic differential equation with no drift.
Here we have adopted the convention that $\widehat{u^j}$ means the omission
of the vector $u^j$ in the tensor.
For a general $q$-vector $U$ in $\wedge^qT_{x_0}M$ this could be written
$$div\left(\overline{\wedge^qT_{x}\xi}_\sigma(U)\right)_{t,\dots, t}
=-{1\over t}  \overline{\wedge^{q-1}(T\xi_t)
\iota_{\int_0^t \langle T\xi_s -, dx_s\rangle}U}.$$
 Equation (\ref{1}) will follow, as in \cite{Elworthy-Li-forms-CR},
by calculating the conditional expectations. Note that we have used
another convention: if $\ell$ is a linear map on a vector space
 $$\iota_\ell(b^1\wedge\dots\wedge b^q)
 =\sum_{j=1}^q (-1)^{j+1}\ell(b^j)b^1
   \wedge\dots\wedge\widehat{b^j}\wedge\dots \wedge b^q.$$

In \cite{Elworthy-Li-Hodge-2}, see also \cite{Elworthy-Li-Hodge-1},
a general theory for forms on path spaces is given. There we define
Hilbert spaces $\HH^{(q)}$, densely included in
$\wedge^q T_\sigma C_{x_0}M$ and show that suitably regular sections of
$\HH^{(q)}$ do have a divergence and give evidence to suggest that
the divergence is a section of $\HH^{(q-1)}$, with proofs for $q=2,3$. This
is the major step in the construction of an $L^2$ deRham complex on
$C_{x_0}M$, with an $L^2$ Hodge decomposition and analogue of the
 Hodge-Kodaira Laplacian. For this construction we take a stochastic
differential equation (\ref{5}) on $M$ whose solutions are Brownian
motions and with its associated connection $\breve \nabla$ the Levi-Civita
connection, for example by taking an isometric embedding in $\R^m$
and  $A\equiv 0$. We take $(B_\cdot)$ to be the canonical Brownian
motion defined up to time $T$ so that  $C_0\R^m$ with Wiener measure $\PP$
is the underlying probability space. For fixed $x_0\in M$ there is the
It\^o map
\begin{eqnarray*}
\I: C_0\R^m &\to& C_{x_0}M\\
\I(\omega)_t&=&\xi_t(x_0,\omega).
\end{eqnarray*}
It is smooth in the sense of Malliavin Calculus with $H$-derivative
$$T_\omega\I: H\to TC_{x_0}M$$
defined for almost all $\omega\in C_0\R^m$.  Here $H$ is the Cameron-Martin
space; $H=L_0^{2,1}\left([0,T];\R^m\right)$. We have the formula
of Bismut
\begin{equation}
\label{formula-TI}
T\I(h)_t=T\xi_t\int_0^t \left(T\xi_s\right)^{-1}X(x_s)\dot h_s\; ds.
\end{equation}
By an $L^2$ $H$-$q$-vector field on $C_0\R^m$ we mean an element
$U\in L^2\left(C_0\R^m;\wedge^qH\right)$ where $\wedge^q H$ refers to
the completed tensor product using the usual Hilbert space cross norm.
Given such a vector field $U$ we can define
$$(\wedge^qT\I)(U): C_0\R^m \to \wedge^q TC_{x_0}M.$$
At $\omega\in C_0\R^m$ this is in $\wedge^q T_{x_\cdot(\omega)}C_{x_0}M$.

To obtain a $q$-vector field on $C_{x_0}M$ define, for $\sigma\in C_{x_0}M$
$$\overline{(\wedge^qT\I)(U)}_\sigma
=\E\left\{(\wedge^qT\I)(U)\left| x_\cdot=\sigma\right.\right\}.$$
Taking constant $U$, that is $U(\omega)\equiv h$, some $h\in H$ we obtain
continuous linear maps
$$\overline{(\wedge^qT\I)}_\sigma: \wedge^qH\to \wedge^qT_\sigma C_{x_0}M$$
defined by $$\overline{(\wedge^qT\I)}_\sigma(h)
 =\overline{(\wedge^qT\I)(h)}_\sigma.$$
We can then define $\HH^{q}_\sigma=Image(\overline{\wedge^q T\I}_\sigma)$
with its induced Hilbert space structure, inner product
 $\langle, \rangle_\sigma$, so that $(\overline{\wedge^q T\I})_\sigma$
becomes an isometry. These spaces are defined for almost all
 $\sigma\in C_{x_0}M$. It is shown in \cite{Elworthy-Li-Hodge-2} that
 they depend only on the Riemannian structure of $M$. The fact that
suitably regular $L^2$ sections of $\HH^{q}$ have a divergence can be
 deduced from the theory for Wiener spaces given by Shigekawa
\cite{Shigekawa-Hodge}.

The following is a specialization of his basic result to primitive
adapted $q$-vector fields.
\begin{theorem} [Shigekawa \cite{Shigekawa-Hodge}]
\label{theorem-Shigekawa}
For $q>1$, let $\phi:C_0\R^m\to (\wedge^q H)^*$ be an
H-differential $q$-form on Wiener space and $h^i: C_0\R^m \to H$,
$i=1$ to $q$ adapted vector fields such that
$h:=h^1(\cdot)\wedge\dots \wedge h^{q+1}(\cdot)$ is in
$\D^{1,2}(\wedge^{q+1}H)$. If also $\phi\in
\D^{1,2}\left((\wedge^qH)^*\right)$ then
$$\int_{C_0\R^m} d\phi(h)\; d\PP
=-\int_{C_0\R^m} \phi\left(div h\right)\, d\PP$$ where
\begin{eqnarray*}
div h&=&\sum_{j=1}^q (-1)^{j} \int_0^T \langle \dot h^j_s,
dB_s\rangle \left(h^1\wedge \dots \wedge \widehat{h^j}\wedge \dots
\wedge h^q \right)\\
&&-\sum_{1\le i<j\le q+1}(-1)^{i+j}[h^i,h^j]\wedge h^1 \wedge\dots
\wedge \widehat {h^i}\wedge\dots \wedge \widehat {h^j}\wedge\dots
\wedge  {h^{q+1}}.
\end{eqnarray*}
\end{theorem}
The assumptions of primitivity of $h$ or adaptness of the components could
have been omitted with Skorohod integrals replacing It\^o integrals,
 see \cite{Shigekawa-Hodge}.
\bigskip

The extension to forms on $C_{x_0}M$ in
\cite{Elworthy-Li-Hodge-1}, \cite{Elworthy-Li-Hodge-2} and the
proof we give below for formulae (1) and (2) is based on the
following observation:

\begin{theorem}\label{theorem-2.2}
Let $\psi$ be a smooth cylindrical $q$-form on $C_{x_0}M$. Then
for $h\in \D^{1,2}(\wedge^{q+1}H)$ of the form $h=h^1\wedge \dots
\wedge h^{q+1}$ with each $h^i$ adapted, the $q$+$1$-vector field
$\overline{\wedge^{q+1}T\I(h)}$ on $C_{x_0}M$ has a divergence:
\begin{equation}
\int_{C_{x_0}M} d\psi\left(\overline{\wedge^{q+1}T\I(h)}\right)
d\mu_{x_0} =-\int_{C_{x_0}M} \psi\left(div
\overline{\wedge^{q+1}T\I(h)}\right) d\mu_{x_0}
\end{equation}
and
\begin{equation}
\label{18}
div \left(\overline{\wedge^{q+1}T\I(h)}\right)
= \overline{\wedge^{q}T\I(div h)}.
\end{equation}
\end{theorem}

{\it Proof.}
  The pull back form $\I_*(\psi)$ is in
  $\D^{1,2}\left((\wedge^qH)^*\right)$. The integration by parts
  formula of Shigekawa, Theorem \ref{theorem-Shigekawa},  applies to give:
  \begin{eqnarray*}
&&\int_{C_{x_0}M} d\psi\left(\overline{\wedge^{q+1}T\I(h)}\right)
d\mu_{x_0}
=\int_{C_0\R^m} d(\I_*\psi)(h)\; d\PP\\
&=&-\int_{C_0\R^m} \I_*\psi\left(div h\right)\, d\PP
=-\int_{C_{x_0}M}\psi\left(\overline{\wedge^{q}T\I(div h)}\right)
d\mu_{x_0}. \\
\end{eqnarray*}
\hfill\rule{3mm}{3mm}

For other discussions of differential forms
 on path spaces see \cite{Jones-Leandre91}, \cite{Kusuoka-forms}
 and \cite{leandre-survey} together with the references there.

\section{Some families of $H$-valued vector fields}

{\bf 3.1 \hskip 5pt The Lie brackets.}
 Below we identify a family of vector fields whose Lie
 brackets behave particularly well.

\begin{proposition}\label{proposition-3.1}
Given tangent vectors $b_i \in T_{x_0}M, i=1,2$ and   $\rho: [0,T]\to \R$,
 a non-random $L^1$ function, set:
\begin{equation}\label{9}
h_{\cdot }^i=\int_0^{\cdot } \rho(r)Y_{x_r}\left( T\xi _r(b^i)\right) \;dr.
\end{equation}
Then the Lie bracket of the two vector fields $h^i$ on the Wiener space
 $C_0\left(\R^m\right)$ is given by:
\begin{equation}
\begin{split}
\frac d{dt}\left[ h^1,h^2\right]_t
 &=\rho(t) \left(\int_0^t \rho(r)\, dr \right)
         (dY)\left(T\xi _t(b^1), T\xi_t (b^2)\right) \\
 &+\rho(t)Y_{x_t}\left( T\xi_t \int_0^t  \rho(r) (T\xi_r)^{-1}  X(x_r)dY
    \left(T\xi_r(b^1),T\xi _r(b^2)\right)\, dr\right). \end{split}
\end{equation}
and so
\begin{equation}
\label{20}
T\I_t([h^1,h^2])
=\left(\int_0^t\rho(s)ds\right) T\xi_t\left(\int_0^t\rho(s) T\xi_s^{-1}
X(x_s)dY\left(T\xi_s(b^1), T\xi_s(b^2)\right)ds\right).
\end{equation}
\end{proposition}

{\it Proof.}  Denote by $D$ differentiation on Euclidean spaces, $D^H$ the
H-derivatives of Wiener functionals,  $\nabla$ covariant differentiation
using the Levi-Civita connection on $M$, and   $\nabla ^H$ the corresponding
covariant H-derivatives. Below we regularly abuse notation by omitting
$\omega$ in the argument, especially with $\xi_t$ and $T\xi_t$.
First observe that
\begin{equation}
\label{22a}
\begin{split}
T{\I}_t\left(h^i\right) (\omega) =&T\xi _t\int_0^t
(T\xi_r)^{-1} X(x_r(\omega))(\dot h_r^i)\,dr \\
 =&\left(\int_0^t \rho(r)dr\right)T\xi_t(b^i)
\end{split}
 \end{equation} and that
\begin{equation}
\label{derivative-of-Lie-bracket}
\begin{split} \frac d{dt}\left[ h^1,h^2\right] _t(\omega )
&=\frac d{dt}\left(
D^Hh_t^2(\omega )(h^1(\omega ))-D^Hh_t^1(\omega )(h^2(\omega ))\right) \\
&=D^H\dot{h}_t^2(\omega )(h^1(\omega ))-D^H\dot{h}_t^1(\omega )
(h^2(\omega)).
\end{split}
\end{equation}
Now
\begin{eqnarray*}
&&D^H\dot{h}_t^2(\omega )(h^1(\omega )) \\
&=&D^H\left( Y_{x_t}\left( T\xi _t(\rho(t)b^2)\right)
\right)(\omega) (h^1(\omega )) \\
&=&\rho(t) \left(\nabla_{h^1(\omega )}^H Y_{x_t}\right) \left(T\xi_t(b^2)\right)
  +   \rho(t) Y_{x_t(\omega )}\left(\nabla^H _{h^1(\omega )}
    \left(T\xi _t(b^2)\right) \right) \\
&=&\rho(t) \left(\nabla_{T_\omega{\I}_t(h^1)} Y\right)
  \left(T\xi _t(b^2)\right) +\rho(t)  Y_{x_t(\omega )}
\left( \nabla^H_{h^1(\omega )}\left(T\xi _t(b^2)\right)\right).
\end{eqnarray*}
By the expression (\ref{22a}) for $T\I_t(h^1)$ , we have
\begin{equation}
\label{derivative-of-h}
\begin{split}
D^H\dot{h}_t^2(\omega )(h^1(\omega ))
=&\rho(t)\left(\int_0^t\rho(r)dr\right)
   \nabla_{T\xi _t(b^1)} Y\left(T\xi_t(b^2)\right)\\
 & +\rho(t)  Y_{x_t(\omega )}\left(\nabla _{h^1(\omega )}^H\left(T\xi_t(b^2)\right)\right).
\end{split}
\end{equation}
To calculate the second term, consider the map
\begin{eqnarray*}
&\tilde {\I}_t:&M\times C_0\left( R^m\right) \rightarrow M\\
& \tilde \I_t(x,\omega)&=\hskip 8pt\xi_t(x)(\omega)
\end{eqnarray*}
with $T_{(x,\omega)}{\tilde\I_t}:T_xM\times L_0^{2,1}(\R^m)
          \to T_{\xi_t(x)(\omega)}M$
 its derivative in the sense of Malliavin calculus. Then
 $$T_{(x,\omega)} \tilde{\I}_t(u,h)
 =T_x\xi _t \int_0^t\left(T_x\xi_r\right)^{-1}
   X(\xi _r(x)) (\dot h_r)\,dr  +  T\xi_t(u),$$
  and $T_{( x_0,\omega)}{\tilde{\I}_t}(u,h(\omega))=T\I_t(h(\omega))+T\xi(u)$.
 Since $b^2$ is non-random,
\begin{eqnarray*}
&&Y_{x_t(\omega )}\left( \nabla _{h^1(\omega )}^H
   \left(T\xi _t(b^2)\right) \right)
=  Y_{x_t(\omega )}\left(\nabla^H_{h^1(\omega)} (T\xi_t)\right) (b^2)  \\
&=& Y_{x_t(\omega ))} \nabla _{b^2}
    \left[ T_{(-,\omega)}{\tilde{\I}}\left(0,h^1(\omega)\right) \right]_t.
\end{eqnarray*}
But
\begin{eqnarray*}
&& \nabla _{b^2} \left[ T_{(-,\omega)}{\tilde{\I}}
 \left(0,h^1(\omega)\right) \right]_t\\
&=&\nabla_{b^2}\left[ T_{-}\xi _t \int_0^t
   \left(T_{-}\xi_r\right) ^{-1}X(\xi _r(-,\omega ))
      \left( \dot{h}_r^1\left(\omega \right) \right) \,dr \right]  \\
&=&\nabla _{b^2}\left[ T_{-}\xi _t \int_0^t
   \left( T_{-}\xi_r\right) ^{-1} X(\xi _r(-))
   \rho(r)  Y_{x_r}\left(T\xi _r(b^1)\right)
     \,dr \right].
\end{eqnarray*}
Since $\left( T_{-}\xi _r\right) ^{-1}X(\xi _r(-,\omega ))
        Y_{\xi _r(-)}\left(T\_\xi _r\right)=Id$,

\begin{eqnarray*}
&&\nabla _{b^2}\left[T_{\left(-,\omega \right) }{\tilde{\I}}
\left(0, h^1(\omega )\right) \right] _t\\
&=&\left(\int_0^t\rho(r)dr \right) \nabla _{b^2}\left( T\xi _t\right)(b^1)\\
&&-T_{x_0}\xi_t\int_0^t\left(T_{x_0}\xi _r\right)^{-1}X(x_r)
    \nabla _{b^2}\left(\rho(r) Y_{\xi _r(- )}
     \left( T_{-}\xi _r(b^1)\right) \right)  \\
&=& \left(\int_0^t\rho(r)dr \right)
     \nabla _{b^2}\left( T\xi _t\right)(b^1)
     -T_{x_0}\xi _t\int_0^t\rho(r)
   \left(T_{x_0}\xi _r\right)^{-1}
  \nabla_{b^2} \left(T_{-}\xi _r\right)(b^1)\\
&& -T_{x_0}\xi _t\int_0^t \rho(r)\left( T_{x_0}\xi _r\right) ^{-1}
     X(x_r)\nabla _{b^2} \left( Y_{\xi _r(- )}\right)
  \left(T_{x_0}\xi _r(b^1)\right)   \\
  &=&  \left(\int_0^t\rho(r)dr \right)
  \, \nabla _{b^2}\left( T\xi _t\right) \left( b^1\right)
   - T\xi _t\int_0^t \rho(r) \left(T\xi _r\right)^{-1}
  \nabla_{b^2}\left(T\xi _r\right)(b^1) \;dr\\
&& -  T\xi_t \int_0^t \rho(r) \left( T\xi_r\right)^{-1}
     X(x_r)\left(\nabla_{T\xi _r(b^2)} Y\right)
     \left( T\xi _r(b^1)\right) \;dr,
\end{eqnarray*}
 and so
\begin{eqnarray*}
&&\nabla _{b^2}\left[ T_{(-,\omega)}\tilde\I
   \left(0, h^1(\omega )\right) \right] _t-\nabla _{b^1}
  \left[ T_{\left(-,\omega \right) }\tilde{\I}
         \left(0, h^2(\omega )\right)\right] _t \\
&=&\left(\int_0^t \rho(r)dr\right)
 \left( \nabla _{b^2}\left( T\xi _t\right) \left( b^1\right)
-\nabla_{b^1}\left(T\xi _t\right) \left( b^2\right) \right)\\
&& - T\xi _t\int_0^t \rho(r) \left( T\xi _r\right)^{-1}
      \left( \nabla_{b^2}\left( T\xi _r\right)(b^1)
       -\nabla _{b^1}\left( T\xi_r\right)(b^2) \right) \;dr \\
&& -T\xi _t\int_0^t  \rho(r)\left( T\xi _r\right) ^{-1}
   X(x_r) dY\left( T\xi _r(b^2), T\xi _r(b^1)\right)\;dr \\
&=& -T\xi _t  \int_0^t\rho(r) \left(T\xi _r\right)^{-1}
   X(x_r) dY\left( T\xi _r(b^2),T\xi _r(b^1)\right)\;dr.
\end{eqnarray*}
Finally combining (\ref{derivative-of-Lie-bracket}), (\ref{derivative-of-h})
and the last identity, we have

\begin{eqnarray*}
&&\frac d{dt}\left[ h^1,h^2\right]_t(\omega ) \\
&=& \rho(t)\left(\int_0^t\rho(r)dr\right)
    \left(\nabla_{T\xi _t(b^1)} Y\right)
      \left( T\xi _t(b^2)\right)\\
&&-\rho(t)\left(\int_0^t\rho(r)dr\right)
    \left(\nabla_{T\xi _t(b^2)} Y \right)
      \left( T\xi _t(b^1)\right)\\
 &&  -  \rho(t)Y_{x_t}\left( T\xi _t  \int_0^t\rho(r) \left(T\xi _r\right)^{-1}
   X(x_r) dY\left( T\xi _r(b^2),T\xi _r(b^1)\right)\;dr\right)\\
&=&\rho(t)\left(\int_0^t\rho(r)dr\right)
    dY\left(T\xi _t(b^1),T\xi _t(b^2)\right)\\
 &&  + \rho(t)Y_{x_t}\left( T\xi _t  \int_0^t\rho(r) \left(T\xi _r\right)^{-1}
   X(x_r) dY\left( T\xi _r(b^1),T\xi _r(b^2)\right)\;dr\right).
\end{eqnarray*}

To prove (\ref{20}), simply do an integration by parts:

\begin{eqnarray*}
T\I_t([h^1,h^2])&=&T\xi_t\int_0^t\rho(s)\left(\int_0^s\rho(r)dr\right)
(T\xi_s)^{-1}XdY\left(T\xi _s(b^1), T\xi_s (b^2)\right) ds\\
&&+T\xi_t\int_0^t \rho(s)\int_0^s  \rho(r) (T\xi_r)^{-1} XdY
    \left(T\xi_r(b^1),T\xi _r(b^2)\right)\, dr ds\\
    &=&T\xi_t\int_0^t\rho(s)\left(\int_0^s\rho(r)dr\right)
(T\xi_s)^{-1}XdY\left(T\xi _s(b^1), T\xi_s (b^2)\right) ds\\
&&+T\xi_t\left(\int_0^t \rho(s)ds\right) \int_0^t  \rho(r) (T\xi_r)^{-1} XdY
    \left(T\xi_r(b^1),T\xi _r(b^2)\right)\, dr ds\\
    &&-T\xi_t\int_0^t\left( \int_0^s \rho(r)dr\right)  \rho(s) (T\xi_s)^{-1} XdY
    \left(T\xi_s(b^1),T\xi _s(b^2)\right)\,ds\\
&=&\left(\int_0^t\rho(s)ds\right) T\xi_t\left(\int_0^t\rho(s)(T\xi_s)^{-1}XdY
\left(T\xi_s(b^1), T\xi_s(b^2)\right)ds\right).
\end{eqnarray*}

\hfill\rule{2mm}{2mm}

Note that if $\breve T$ is the torsion of the connection $\breve \nabla$,
then  $X(dY)={\breve T}$ by  (\ref{8}). Consequently,
\begin{corollary} For $h^i$ defined by (\ref{8}) in Proposition
  \ref{proposition-3.1},
 \begin{equation}
 \begin{split}X(x_t){d\over dt}[h^1,h^2]_t
 =&\rho(t) \left(\int_0^t \rho(r)\, dr \right)
        \breve T\left(T\xi _t(b^1), T\xi_t (b^2)\right) \\
 &+\rho(t)  T\xi_t \int_0^t  \rho(r) (T\xi_r)^{-1} \breve T
    \left(T\xi_r(b^1),T\xi _r(b^2)\right)\, dr \end{split}
\end{equation}
and
\begin{equation}
\label{24}
T\I_t([h^1,h^2])
=\left(\int_0^t\rho(s)ds\right) T\xi_t\left(\int_0^t\rho(s) T\xi_s^{-1}
\breve T\left(T\xi_s(b^1), T\xi_s(b^2)\right)ds\right).
\end{equation}
\end{corollary}

We shall now return to our examples.

\subsubsection*{Example 3.1. Gradient system.}

In Example 1.1, where the stochastic differential equation is a gradient
system, we have $dY\equiv 0$. Consequently
 $\frac d{dt}\left[ h^1,h^2\right]_t(\omega )$ vanishes and so does
$\left[ h^1,h^2\right]_t(\omega )$. Indeed in this case if we set
$$h^v=\int_0^\cdot Y_{x_r}\left(T\xi_r(v)\right)\;dr$$
for $v\in T_{x_0}M$, we have a  commuting family
$\{h^v: v\in T_{x_0}M\}$  of $\HH$-valued vector fields on $C_0\R^m$.
Note that if $\alpha: M\to \R^m$ is the immersion defining our
stochastic differential equation, then $\alpha\circ \xi_t$ is a random
map of $M$ into $\R^m$ and
$\dot h_t^v=d(\alpha\circ \xi_t)(v)=(\xi_t)^*(d\alpha)(v)$.

\subsubsection*{Example 3.2. The Lie group case.}

For $M$ a Lie group with left invariant stochastic differential equation
corresponding to a bi-invariant metric as  discussed in \S\ref{section-1},
if  $v\in \g$ and $x_0=e$ the identity element, we have
\begin{eqnarray*}
h_\cdot^i&:=&\int_0^\cdot \rho(r) Y^L_{x_r}\left(T\xi_r^L(b^i)\right)\, dr\\
&=&\int_0^\cdot \rho(r) TL^{-1}_{x_r}TR_{x_r}(b^i)\, dr\\
&=&\int_0^\cdot \rho(r) ad(x_r^{-1})(b^i)\,dr
\end{eqnarray*}
In this case
$$T\I_t(h^i)=TR_{x_t}\int_0^t ad(x_r)\dot h_r^i \,dr
=TR_{x_t}\left(\int_0^t\rho(r)dr b^i\right).$$

 Let  $T^R$ and $T^L$ be the torsions for the left and right invariant
connections respectively. From above, (\ref{torsion-of-Lie-groups}),
the fact that $T\xi_t^L(u)=TR_{x_t}(u)$, and by the right invariance of $TR$,
 \begin{eqnarray*}
{d\over dt}[h^1,h^2]_t
 &=&-\rho(t) \left(\int_0^t \rho(r)\, dr \right)
       TL^{-1}_{x_t} TR_{x_t} [b^1,b^2] \\
&&- \rho(t) TL_{x_t}^{-1}TR_{x_t}\int_0^t  \rho(r) [b^1,b^2]\, dr\\
&=&-2\rho(t) \left(\int_0^t \rho(r)\, dr \right) ad(x_t^{-1})[b^1,b^2]\\
&=&-{d \over dt}\left(\int_0^t\rho(r)dr\right)^2 ad(x_t^{-1})[b^1,b^2]\\
&=&-ad(x_t^{-1}){d \over dt}
\left[\int_0^t\rho(r)dr b^1, \int_0^t\rho(r)dr b^2\right].
\end{eqnarray*}
Note that
\begin{eqnarray*}
T\I_t([h^1,h^2])
&=&-TR_{x_t}\int_0^t  ad(x_r)({d \over dr}[h^1,h^2]_r)dr\\
&=&-TR_{x_t}\int_0^t 2\rho(r) \left(\int_0^r \rho(s) ds\right) dr[b^1, b^2] \\
&=&-TR_{x_t} [\int_0^t \rho(r)dr b^1, \int_0^t \rho(r)drb^2].
\end{eqnarray*}

\subsubsection*{Example 3.3.  Lie group as a symmetric space.}
We now consider the stochastic differential equation
(\ref{sde-bi-invariant-Lie-group}):
$$dx_t=TR_{x_t}\circ dB_t -TL_{x_t}\circ dB_t^\prime$$
 on the Lie group $G$ with bi-invariant metric. The derivative flow
of the solution $\xi_t$ is given by $T\xi_t=TR_{(g_t^\prime)^{-1}}TL_{g_t}$.
Note that in our notation the map, giving the equation
(\ref{sde-bi-invariant-Lie-group}),
 $X(g):\g\times \g\to T_xM$ is given by
 $X(g)(\alpha_1,\alpha_2)=TR_g(\alpha_1)-TL_g(\alpha_2)$ with inverse
 $Y:TG\to \g\times \g$ given, for $u\in T_gG$, by
 \begin{eqnarray*}
 Y(u)&=&(TR_g^{-1}(u), -TL_g^{-1}(u))\\
 &=&(Y^R u, -Y^Lu).
\end{eqnarray*}
 Thus if $u, v\in T_gG$, for $T^R$ and $T^L$ the torsion tensors as in
 (\ref{torsion-of-Lie-groups}),
 \begin{eqnarray*}
 dY(u,v)&=&\left( TR_g^{-1}T^R(u,v), TL_g^{-1}T^L(u,v) \right)\\
 &=&\left( [Y^R(u), Y^R(v)], \; [Y^L(u), Y^L(v)] \right)
 \end{eqnarray*}
so
\begin{eqnarray*}
{d\over dt}[h^1,h^2]_t
&=&\rho(t)\left(\int_0^t\rho(r)dr\right)
 \left(\;[Y^R \left(TL_{g_t}TR_{g_t^\prime}^{-1}b^1\right),
                Y^R (TL_{g_t}TR_{g_t^\prime}^{-1}b^2)],\right.\\
&&\hskip 100pt \left.[Y^L (TL_{g_t}TR_{g_t^\prime}^{-1}b^2),
    Y^L (TL_{g_t}TR_{g_t^\prime}^{-1}b^1])\;\right).
    \end{eqnarray*}
Now \begin{eqnarray*}
Y^R \left(TL_{g_t}TR_{(g_t^\prime)}^{-1}b^i\right)
&=&\left(TR_{g_t}(g_t^\prime)^{-1}\right)^{-1}TL_{g_t} TR_{g_t^\prime}^{-1}b^i\\
&=&\left( TR_{g_t^\prime}^{-1}\circ TR_{g_t}\right)^{-1}
TL_{g_t} TR_{g_t^\prime}^{-1}b^i\\
&=&ad(g_t)b^i.
\end{eqnarray*}
while
\begin{eqnarray*}
Y^L \left(TL_{g_t}TR_{(g_t^\prime)}^{-1}b^i\right)
&=&\left(TL_{g_t}(g_t^\prime)^{-1}\right)^{-1}
   TL_{g_t} TR_{g_t^\prime}^{-1}b^i\\
 &=&TL_{g_t^\prime}TR_{g_t^\prime}^{-1} b^i=ad(g_t^\prime)b^i,
 \end{eqnarray*}
so
\begin{eqnarray*}
{d\over dt}[h^1,h^2]_t
&=&\rho(t)\left(\int_0^t\rho(r)dr\right)
\left(ad(g_t)[b^1,b^2], ad(g_t^\prime)[b^1,b^2]\right)\\
&\in& \g\times \g.
\end{eqnarray*}
In this case $h_t^i=\left(\int_0^t\rho(r)ad(g_r)b^i\, dr,
-\int_0^t\rho(r)ad(g_r^\prime)b^i\,dr\right)$.
\bigskip

\noindent
{\bf 3.2 \hskip5pt A family of $\HH^2$-valued vector fields. }
In general when the associated connection $\breve \nabla$ for (\ref{5})
is the Levi-Civita connection $\nabla$ and $A\equiv 0$ we expect that
$\overline{\wedge^qT\I(h^1\wedge \dots \wedge h^q)}_\sigma \in \HH^q_\sigma$
for almost all $\sigma \in C_{x_0}M$. It follows from \cite{Elworthy-Li-Hodge-1}
that this holds for $q=1$ and from \cite{Elworthy-Li-Hodge-2} for $q=2$ with sufficiently regular $h$, and a
proof for $q=2$ for our special $h^i$'s  is given below in Proposition
\ref{proposition-3.3} as an illustrative
example. For general $q$ we have,
$$\overline{\wedge^q(T\I)(h^1\wedge \dots \wedge h^q)}
=\E\left\{\wedge^q(T\xi_\cdot)\left(b^{\rho,1}_\cdot\wedge \dots \wedge
b_\cdot^{\rho,q}\right) \left| x_s: 0\le s \le T \right.\right\}$$
where $b_t^{\rho,j}=\left(\int_0^t \rho(r)dr\right) b^j$ and
$(T\xi_\cdot)(b_\cdot^{\rho,j})$ is just the vector field \\
$\{(T\xi_t)(b_t^{\rho,j}): 0\le t \le T\}$
along $\xi_\cdot(x_0)$. From \cite{Elworthy-LeJan-Li-book} or
\cite{Elworthy-Yor}  for gradient systems, we can deduce that for
 $0<t_1<\dots <t_q\le T$
\begin{equation}
\label{25}
\begin{split}
&\overline{\wedge^qT\I(h^1\wedge \dots \wedge h^q)}_{t_1,\dots,t_q}\\
=&\int_0^{t_1}\rho(r)dr\int_0^{t_2}\rho(r)dr\dots \int_0^{t_q}\rho(r)dr
 \cdot \left(\1\otimes \dots \otimes W_{t_q}^{t_{q-1}} \right)\\
&\left(\1\otimes \dots \otimes {W^{(2)}}^{t_{q-2}}_{t_{q-1}}\right)
\dots W_{t_1}^{(q)}(b^1\wedge \dots \wedge b^q).
\end{split}
\end{equation}
Here ${W^{(q)}}^s_t=W_t^{(q)}\circ (W_s^{(q)})^{-1}:
 \wedge^qT_{x_s}M \to \wedge^qT_{x_t}M$
 with $W^s_t={W^{(1)}}^s_t$, for  $W^{(q)}$ the damped parallel transport
 defined by (\ref{3}).

 Let us recall the characterization of $\HH^1$ and $\HH^2$. For $q=2$ it is shown
  in \cite{Elworthy-Li-Hodge-2}, that
 $$\HH^{2}_\sigma=\left\{U+Q_\sigma(U) \quad |\quad
 U\in \wedge^2 \HH^{1}_\sigma\right\}\subset \wedge^2T_\sigma C_{x_0}M $$
 where $Q_\sigma: \wedge^2\HH^1_\sigma \to \wedge^2 T_\sigma C_{x_0}M$
 is defined by
\begin{equation}\label{Q}
Q(U)_{s,t}=(\1\otimes W^s_t) W_s^{(2)} \int_0^s (W_r^{(2)})^{-1}
 \RR(U_{r,r})\, dr, \hskip 15pt 0\le s \le t \le T
\end{equation}
 for $\RR: \wedge^2TM \to \wedge^2 TM$ the curvature operator.

 The space $\HH^1$ is  the `Bismut tangent space', \cite{Jones-Leandre91},
but with a different Hilbert space structure. In fact elements of
$\HH^1$ are exactly vector fields $v$ in $M$ along $\sigma$
 such that $\|v\|_\sigma:=\int_0^T\left \vert {\D\over dt}v_t\right
   \vert^2_{\sigma(t)} \; dt<\infty$ where
\begin{equation}
{\D\over d t}v_t:={D\over d t}v_t+{1\over 2}Ric^{\#}(v_t)
=W_t^{-1}{d\over dt}W_t^{-1}V_t.
\end{equation}

From this it follows that $\wedge^2 \HH^1_\sigma$ consists of those $U$
in $\wedge^2T_\sigma C_{x_0}M$ for which
$$\left(\wedge^2(W_\cdot)^{-1}(U)\right)_{s,t}
=\int_0^s\int_0^t H(r_1,r_2) \, dr_1 dr_2$$
for some $H$ in $L^2\left([0,T]\times [0,T]; T_{x_0}M\otimes T_{x_0}M\right)$.

\begin{proposition}
\label{proposition-3.3}
Let $\lambda\in L_0^{2,1}\left([0,T];\R\right)$.  Take $V$ in
$\wedge^2T_{x_0}M$ and $V^\lambda \in \wedge^2 L_0^{2,1}T_{x_0}M$ defined by
$V^\lambda_{s,t}=\lambda(s)\lambda(t)V$.
If $A\equiv 0$ in (\ref{5}) and the associated connection is the Levi-Civita
connection, the vector field $Z$, where
$$Z=\E\left\{\wedge^2(T\xi_\cdot)(V^\lambda) \vert x_s, 0\le s \le T\right\}$$
forms a section of $\HH^2$. It is equal to $U+Q(U)$ for $Q$ as above
and $U$ given by
\begin{equation}
\label{28}
\begin{split}
U_{s,t}&=(\1\otimes W^s_t)\lambda(s)\lambda(t) W_s^{(2)}(V)\\
&-(W_s\otimes W_t)      \int_0^s    \lambda^2(r) \wedge^2(W_r)^{-1}
\RR\left( W_r^{(2)}(V)\right)\, dr,
\end{split}
\end{equation}
for all $0\le s\le t\le T$. It has divergence
\begin{equation}
\label{divergence}
\begin{split}
(div Z)_t=&-\int_0^T \left(\langle \frac{\D} {\partial r}-, dx_r\rangle
\otimes \1\right) Z_{r,t}\\
&+W_t\int_0^t\left(\langle -, dx_r \rangle
\otimes (W_r)^{-1})\RR(\1\otimes W^r_t)\right)^{-1}Z_{r,t}.
\end{split}
\end{equation}
(For an interpretation of the, apparently adapted, integrals, see the proof below.)
\end{proposition}

{\it Proof.}
From (\ref{25}), for $0\le s\le t\le T$,
\begin{equation}
\label{32a}
Z_{s,t}=\lambda(s)\lambda(t)(\1\otimes W^s_t)W_s^{(2)}(V).
\end{equation}
So
\begin{equation}\label{29}
\left[(\wedge^2 W_\cdot)^{-1} Z
            \right]_{s,t}
=\lambda(s)\lambda(t)(\wedge^2 W_s)^{-1}W_s^{(2)}(V).
\end{equation}
Setting
\begin{equation}\label{30}
k(r)=\wedge^2(W_r^{-1})W_r^{(2)}(V)\in \wedge^2T_{x_0}M
\end{equation}
we see that
$$\left[(\wedge^2 W_\cdot)^{-1} Z
            \right]_{s,t} =\lambda(s)\lambda(t)k(s\wedge t),
\hskip 15pt \hbox{ for all $s,t\in [0,T]$}.$$

Set $$\tilde U=\wedge^2(W^{-1}_\cdot)(U)$$
We only need to show that
\begin{equation}\label{31}
\tilde U+(\wedge^2 W_\cdot^{-1})Q(U) =\lambda(s)\lambda(t)k(s\wedge t)
\end{equation}
 and that $U$ belongs to $\wedge^2 \HH^1$.

  Now for any  $K \colon [0,T]\times \Omega\to \wedge^2TM$ over
 $\{x_t: 0\le t\le T\}$
 \begin{equation}\label{32}
 \wedge^2(W_r){d \over dr}\left(\wedge^2(W_r)^{-1}K_{r}\right)
 =\RR (K_{r})+  W^{(2)}_r{d \over dr}[(W_r^{(2)})^{-1}K_r]
 \end{equation}
when the derivative exists, so
$k^\prime(s)=\wedge^2(W_s^{-1})\RR \left(W_s^{(2)}(V)\right)$
and consequently by (\ref{28})
\begin{equation}\label{33}
\tilde U_{s,t}=\lambda(s)\lambda(t)k(s\wedge t)
-\int_0^{s\wedge t} \lambda^2(r)k^\prime(r)\, dr,
\end{equation}
 for $s,t\in [0,T]$. Integrating by parts,
\begin{eqnarray}
\nonumber
\tilde U_{s,t}&=&\lambda(s)\lambda(t)k(s\wedge t)
  -\lambda^2(s\wedge t)k(s\wedge t)
 + 2 \int_0^{s\wedge t}\lambda^\prime(r)\lambda(r)k(r)\,dr\\
&=&2 \int_0^{s\wedge t}\lambda^\prime(r)\lambda(r)k(r)\,dr
  +(\lambda(s)\lambda(t)-\lambda^2(s\wedge t))k(s\wedge t).
  \label{36}
\end{eqnarray}
From this it is  easy to see that
 \begin{equation}
\label{34}
\tilde U_{\bar s \bar t}=\int_0^{\bar s}ds \int_0^{\bar t}dt
\left(\lambda^\prime(s)\lambda^\prime(t)k(s\wedge t)
+\lambda(s\wedge t)\lambda^\prime(s\vee t) k^\prime(s\wedge t)\right)
\end{equation}
since for  $\bar s \le \bar t$ one has
\begin{equation}
\begin{split}
&\int_0^{\bar s}ds \int_0^{\bar t}dt
\left(\lambda^\prime(s)\lambda^\prime(t)k(s\wedge t)
+\lambda(s\wedge t)\lambda^\prime(s\vee t) k^\prime(s\wedge t)\right)\\
=&\int_0^{\bar s}\left\{ \int_0^{s}
\left(\lambda^\prime(s)\lambda^\prime(t)k(t)
+\lambda(t)\lambda^\prime(s) k^\prime(t)\right)dt\right.\\
& \left.+\int_s^{\bar t}  \left(\lambda^\prime(s)\lambda^\prime(t)k(s)
+\lambda(s)\lambda^\prime(t) k^\prime(s)\right)dt\right\}\\
=&   \int_0^{\bar s}\lambda^\prime(s)\lambda(s)k(s)\,ds
+\int_0^{\bar s}(\lambda(\bar t)-\lambda(s)){d \over ds}(\lambda(s)k(s))\,ds\\
=& 2 \int_0^{\bar s}\lambda^\prime(s)\lambda(s)k(s)\,ds
+(\lambda(\bar t)-\lambda(\bar s))\lambda(\bar s)k(\bar s).
\end{split}
\end{equation}
This shows that $U\in \wedge^2\HH^1$. Next note that
$$\tilde U_{s,s}=2\int_0^s \lambda^\prime(r)\lambda(r)k(r)\, dr$$
by (\ref{36}). By the definition,  (\ref{Q}),  of $Q$, (\ref{30}),
and using (\ref{32}) again  we see that  for $0\le s\le t\le T$,

\begin{eqnarray*}
&&\left[\wedge^2(W_\cdot^{-1})Q(U)\right]_{s,t}
:=\wedge^2(W_s^{-1})W_s^{(2)}\int_0^s (W_r^{(2)})^{-1}
\RR\left(U_{rr}\right) dr\\
&=&\wedge^2(W_s^{-1})W_s^{(2)}\int_0^s (W_r^{(2)})^{-1}\wedge^2(W_r)
 {d\over dr}\tilde U_{r,r}\, dr\\
 &&-\wedge^2(W_s^{-1})W_s^{(2)}\int_0^s{d \over dr}
 \left[(W_r^{(2)})^{-1} U_{rr}\right]\,dr\\
 &=&\wedge^2(W_s^{-1})W_s^{(2)}\int_0^s 2 (W_r^{(2)})^{-1}\wedge^2(W_r)
  \lambda^\prime(r)\lambda(r)k(r)\, dr-\tilde U_{s,s}\\
&=&\wedge^2(W_s^{-1})W_s^{(2)}\int_0^s 2
  \lambda^\prime(r)\lambda(r)V\, dr-\tilde U_{s,s}\\
  &=&\lambda^2(s)k(s)-2\int_0^s \lambda^\prime(r)\lambda(r)k(r)\, dr \\
&=&\int_0^s  [\lambda(r)]^2 k^\prime(r) \, dr.
\end{eqnarray*}
Finally add (\ref{33}) to the above to obtain the required identity
(\ref{31}).

To obtain (\ref{divergence}) we can assume $V$ to be primitive,
$V=b^1\wedge b^2$ say. By (\ref{18}) the divergence of $Z$ is just the
conditional expectation $\bar J$ of $J$ for $J=T\I(div(h^1\wedge h^2))$
where $\dot h_t^j=\dot \lambda(t)Y_{x_t}T\xi_t(b^j)$.
By Theorem \ref{theorem-Shigekawa} and Proposition \ref{proposition-3.1}
\begin{eqnarray*}
J_t&=&T\I_t\left(
-\int_0^T\langle \dot h_r^1, dB_r\rangle h_\cdot^2
+\int_0^T\langle \dot h_r^2, dB_r \rangle h_\cdot^1\right)\\
&=&J_t^0+J_t^1
\end{eqnarray*}
where
\begin{eqnarray*}
J_t^0&=&T\I_t\left(-\int_t^T\langle \dot h_r^1, dB_r\rangle h^2_\cdot
+\int_t^T\langle \dot h_r^2, dB_r\rangle h_\cdot^1\right)\\
&=&-\int_t^T\left\langle \dot \lambda(r)T\xi_r^t(T\xi_t(b^1)), dx_r
\right\rangle  \lambda(t) T\xi_t(b^2)\\
&&+\int_t^T \left\langle \dot \lambda(r) T\xi_r^t (T\xi_t(b^2)), dx_r
 \right\rangle \lambda(t)T\xi_t(b^1).
\end{eqnarray*}
Therefore
$$\E\left\{ J_t^0\vert \; \F^{x_0}\vee \F_t\right\}
=-\lambda(t)\int_t^T \dot\lambda(r)\iota_{\langle W_r^t-,dx_r\rangle}
\wedge^2(T\xi_t)(b^1\wedge b^2)$$
giving
$$\bar J_t^0=\int_t^T\left(\1\otimes \left\langle \frac{\D}{\partial r}-,
dx_r\right\rangle \right) \, Z_{t,r}$$
by (\ref{32a}).

The second term $J_t^1$ could be treated as in Corollary \ref{cor-invariant}
below, but to get the form (\ref{divergence}) observe
$$J_t^1=- \int_0^t \lambda(t)\dot \lambda(r)
\left(\langle-, dx_r\rangle \otimes \1\right)
\left(\wedge^2(T\xi)(b^1\wedge b^2)\right)_{rt}$$
and so by (\ref{25})
$$\bar J_t^1=-\int_0^t \lambda(t)\dot \lambda(r)
\left(\langle-, dx_r\rangle \otimes \1\right)
\left(\1\otimes W_t^r)W_r^{(2)}(b^1\wedge b^2)\right)_{rt}.$$
Here we can interpret the, apparently adapted, integral by using
$W^r_t=W_t(W_r)^{-1}$ and taking the term $(\1\otimes W_t)$
outside of the integral.
Now
\begin{eqnarray*}
(\frac{\D}{\partial r}\otimes \1)Z_{rt}
&=&\frac{D}{\partial r}\left[ \lambda(r)\lambda(t)
(\1\otimes W^r_t)W^{(2)}_r(b^1\wedge b^2)\right]\\
&&+\frac{1}{2}(Ric^{\#}\otimes \1)Z_{rt}
\end{eqnarray*}
giving
\begin{eqnarray*}
&&\dot\lambda(r)\lambda(t)(\1\otimes W^r_t)W^{(2)}_r(b^1\wedge b^2)\\
&=&(\frac{\D}{\partial r}\otimes \1)Z_{rt}-\frac{1}{2}(Ric^{\#}\otimes \1)Z_{rt}\\
&&-\lambda(r)\lambda(t)(1\otimes W_t)\frac{D}{\partial r}
(\1\otimes (W_r)^{-1})W_r^{(2)}(b^1\wedge b^2),
\end{eqnarray*}
but
\begin{eqnarray*}
&&\frac{D}{\partial r}(\1\otimes (W_r)^{-1})W_r^{(2)}(b^1\wedge b^2)\\
&&=(W_r\otimes \1)\frac{D}{\partial r}(\wedge^2 W_r^{-1})
W_r^{(2)}(b^1\wedge b^2)
-\frac{1}{2}(Ric^{\#}\otimes \1)(\1\otimes W_r^{-1})
W_r^{(2)}(b^1\wedge b^2)\\
&&=(\1\otimes W_r^{-1})\RR W_r^{(2)}(b^1\wedge b^2)
-\frac{1}{2}(Ric^{\#}\otimes \1)(\1\otimes W_r^{-1})W_r^{(2)}(b^1\wedge b^2)
\end{eqnarray*}
by (\ref{32}), and so
$$\dot\lambda(r)\lambda(t)(\1\otimes W^r_t)W^{(2)}_r(b^1\wedge b^2)
=(\frac{\D}{\partial r}\otimes \1)Z_{rt}
-(\1\otimes W_t^r)\RR (\1\otimes W_t^r)^{-1}Z_{rt}.$$
The result follows noting that $(\1\otimes W^r_t)^{-1}Z_{rt}$
is $\F_r^{x_0}$-measurable and the factor $(\1\otimes W_t)$
in $(\1\otimes W^r_t)=(\1\otimes W_t)(\1\otimes (W_r)^{-1})$
can be taken out of the stochastic integrals to make them It\^o
integrals.

\hfill \rule{3mm}{3mm}

\bigskip

\noindent{\bf Remark 3.3.}
\vskip 2pt

(1).   The formula (\ref{28}) in Proposition \ref{proposition-3.3} demonstrates the
alternative characterisation of $\HH^2$ from \cite{Elworthy-Li-Hodge-2}
 which states that $V$ is in $\HH^2$ if and only if $V-\pathR(V)$ is in
 $\wedge^2 \HH^1$  where $\pathR$ is the curvature operator of the damped
 Markovian connection on $\HH^1$.
\bigskip

(2). Note that for the proof of (\ref{28}) we could allow
$V: C_{x_0}M\to \wedge^2T_{x_0}M$ to be nonconstant provided it is in
$L^{1+\epsilon}$ for some $\epsilon>0$. Our filtration in calculating
(\ref{25}) can be chosen to be $\{\F_t \vee \F^{x_0}: 0\le t\le T\}$ as
 in \cite{Elworthy-Li-Ito}.
Now let $V: [0,T]\times [0,T] \to \wedge^2 T_{x_0}M$ be such that
$V_{s,t}=\int_0^s\int_0^t Z_{a,b}\, da\, db$ some
$Z\in L^2\left([0,T]\times [0,T]\times C_{x_0}M\to \wedge^2T_{x_0}M\right)$
symmetric in $(s,t)\in [0,T]\times [0,T]$. As a corollary we see,
 for $0\le s\le t \le T$,
$$\E\{\wedge^2(T\xi_\cdot)(V)\vert \,x_s: 0\le s\le t\le T\}=U+Q(U)$$
for $$U_{s,t}=(\1\otimes W^s_t)W^{(2)}_s(V_{s,t})
-(W_s \otimes W_t)
\int_0^s \wedge^2(W_r)^{-1}\RR \left(W^{(2)}_r(V_{rr})\right)\, dr.$$
Moreover $U$ is a section of $\wedge^2\HH^1$ and so
$U+Q(U)$ is a section of $\HH^2$. This follows from the proposition
by polarization for
$V_{s,t}=\left(\lambda^1(s)\lambda^2(t)+\lambda^2(s)\lambda^1(t)\right)b$,
 some $b\in L^2(C_{x_0}M;\wedge^2T_{x_0}M)$ and then for general $V$ by
 continuity since we can consider such $V$ as elements of the completed
tensor product
$[\odot^2L_0^{2,1}([0,T];\R)]\otimes L^2(C_{x_0}M;\wedge^2T_{x_0}M)$,
where $\odot$ refers to the symmetric tensor product. Note the mapping
\begin{eqnarray*}
&&\odot^2L_0^{2,1}\left([0,T];\R\right)
 \otimes L^2(C_{x_0}M;\wedge^2T_{x_0}M)\\
&&\stackrel{\Theta}{\longrightarrow}
L^2\left(C_{x_0}M; \wedge^2(L_0^{2,1}([0,T];\R)\otimes T_{x_0}M)\right)\\
&&\simeq L^2\left(C_{x_0}M; \wedge^2(L_0^{2,1}(T_{x_0}M)\right)
\end{eqnarray*}
given by
\begin{eqnarray*}
&&\Theta\left((u\odot v)\otimes f(\alpha\wedge\beta)\right)(\sigma)\\
&&={1\over 2} f(\sigma)\left[(u\otimes \alpha)\wedge (v\otimes \beta)
+(v\otimes \alpha)\wedge (u\otimes \beta)\right],\sigma\in C_{x_0}M,
\end{eqnarray*}
 where $f(\cdot)(\alpha\wedge \beta)
\in L^2\left(C_{x_0}M;\wedge^2T_{x_0}M\right)
\simeq L^2(C_{x_0}M;\R)\otimes \wedge^2 T_{x_0}M$,\\
$\alpha, \beta\in T_{x_0}M$, $f\in L^2(C_{x_0}M;\R)$.

In this generality the stochastic integrals in (\ref{divergence})
will need to be treated more carefully and if $V$ is non-random an
additional term involving its H-derivative will be involved in the divergence.
\bigskip

\noindent{\bf 3.3 \hskip 5pt
From Integration by parts to  Bismut type formulae.}
First we shall give an extension of the non-intrinsic formula (\ref{4}).
For this consider a general non-degenerate stochastic differential equation
(\ref{5}) with smooth coefficients. For $\{\xi_t: 0\le t <\infty\}$, its
flow of diffeomorphisms of the manifold, consider the semigroup
described in \S\ref{section-1}, given by
$$P_t\phi=\E \xi_t^*(\phi).$$
Since $M$ is compact $d(P_t\phi)=P_t(d\phi)$ when $\phi$ is $C^1$.
The generator $\A$ of $\{P_t\}_{t\ge 0}$ is given on smooth forms by
(\ref{11}) above.

\begin{proposition}
For any bounded measurable form $\phi$ on $M$ and $L^1$ function
 $\rho:[0,t]\to \R$ with $\int_0^t \rho(r)\,dr\not =0$,
\begin{equation}
\label{Bismut-formula}
\begin{split}
&d(P_t\phi)(b^1\wedge \dots \wedge b^{q+1})\\
&=\left(\int_0^t\rho(r)dr\right)^{-1}
  \E  \sum_{j=1}^{q+1}(-1)^{j+1}\int_0^t \rho(s)
\langle T\xi_s(b^j), X(x_s)dB_s\rangle_{x_s} \cdot\\
&   \hskip 130pt \xi_t^*(\phi)(b^1\wedge\dots
 \wedge\widehat{b^j}\wedge \dots \wedge b^{q+1})\\
 &-\left(\int_0^t\rho(r)dr\right)^{-1}
 \sum_{1\le i<j\le q+1}(-1)^{i+j+1}\; \E \xi_t^*(\phi)\\
&  \left(\int_0^t\rho(s)T\xi_s^{-1}
 \left(\breve T\left(T\xi_s(b^i),T\xi_s(b^{j})\right)\right) ds
 \wedge b^1\wedge \dots \wedge\widehat{b^i} \dots
 \wedge\widehat{b^j} \dots \wedge b^{q+1}\right),
\end{split}
\end{equation}
where $\breve T:TM \oplus TM \to TM$ is the torsion of the  connection
$\breve \nabla$ given by (\ref{8}).
\end{proposition}
\noindent
{\it Proof.}
Define $h^i, i=1$ to $q+1$, by (\ref{9}). Set
 $h=h^1\wedge \dots \wedge h^{q+1}$ and $b=b^1\wedge\dots\wedge b^{q+1}$.
Arguing as in Theorem \ref{theorem-2.2}, but without taking conditional
 expectations we see, for $\phi$ a $C^1$ $q$-form,
\begin{eqnarray*}
d(P_t\phi)(b)&=&P_t(d\phi)(b)\\
&=& \E \xi_t^*(d\phi)(b)
=\left(\int_0^T \rho(r)dr\right)^{-q-1} \E
   d\phi\left(\wedge^{q+1}(T\I)(h)_{t,\dots, t}\right)\\
&=& -    \left(\int_0^T \rho(r)dr\right)^{-q-1} \E
  \phi\left(\wedge^{q}(T\I)(div h)_{t,\dots, t}\right).
\end{eqnarray*}
By Shigekawa's result, Theorem \ref{theorem-Shigekawa},
 formula (\ref{Bismut-formula})
follows using Proposition \ref{proposition-3.1},
 equation (\ref{formula-TI}), (\ref{24}), and the fact
that $\int_0^t\langle \dot h^i_s, dB_s\rangle
=\int_0^t\rho(s)\langle T\xi_s(b^i), X(x_s)dB_s\rangle_{x_s}$.
By continuity it also holds for bounded measurable $\phi$.
\hfill\rule{3mm}{3mm}

\bigskip

To obtain an intrinsic formula from (\ref{Bismut-formula}), we shall
take conditional expectations, which can easily be done if
the torsion $\breve T=X(dY)$ is invariant under the flow (e.g.
 for certain homogeneous spaces). In this case
formula (\ref{Bismut-formula}) becomes:

\begin{equation}
\begin{split}
&d(P_t\phi)(b^1\wedge \dots \wedge b^{q+1})\\
&=\left(\int_0^t\rho(r)dr\right)^{-1}
  \E  \sum_{j=1}^{q+1}(-1)^{j+1}\int_0^t \rho(s)
\langle T\xi_s(b^j), X(x_s)dB_s\rangle _{x_s}\cdot\\
&   \hskip 130pt \xi_t^*(\phi)(b^1\wedge\dots
 \wedge\widehat{b^j}\wedge \dots \wedge b^{q+1})\\
 &-\sum_{1\le i<j\le q+1}(-1)^{i+j+1}
  \E\left( \xi_t^*(\phi)\left( \breve T(b^i,b^j)
 \wedge b^1\wedge \dots \wedge\widehat{b^i} \dots
 \wedge\widehat{b^j} \dots \wedge b^{q+1}\right)\right).
\end{split}
\label{40}
\end{equation}
Let $\breve W_t^{A,q}:\wedge^q T_{x_0}M \to \wedge^q T_{x_t}M$,  the damped
 parallel translation of $q$ vectors, be defined by:
\begin{equation}\left\{
\begin{split}
{\hat D\over \partial s}(\breve W_s^{A,q}(V_0))&
 =-{1\over 2}{\breve \RR}_{x_0}^q\left(\breve W_s^{A,q}(V_0)\right)
+d\wedge^q \left(\breve \nabla A\right)\left(\breve W_s^{A,q}(V_0) \right)\\
\breve W_0^{A,q}(V_0)&=V_0.
\end{split}\right.
\end{equation}
Here ${\hat D\over \partial s}$ refers to covariant differentiation
using the connection $\hat \nabla$, the adjoint connection of
$\breve \nabla$. Since the conditional expectation of $\wedge^qT\xi_t$
is given by
$$\overline{\wedge^q T\xi_t}(-)=\breve W_t^{A,q}(-),$$
 by Theorem 3.3.7 in \cite{Elworthy-LeJan-Li-book}, or Theorem A in
 \cite{Elworthy-Yor} for the Levi-Civita connection, we only
 need to worry about the first term on the right hand side of
 equation (\ref{40}).
 Set
 $$U_t=\int_0^t \rho(s) \sum_{j=1}^{q+1}(-1)^{j+1}
 \langle T\xi_s(b^j), X(x_s)dB_s\rangle \cdot \wedge^q T\xi_t
 (b^1\wedge\dots \wedge\widehat{b^j}\wedge \dots \wedge b^{q+1})$$
 and $$V_t^j= \wedge^q T\xi_t
 (b^1\wedge\dots \wedge\widehat{b^j}\wedge \dots \wedge b^{q+1}).$$
 Then,  e.g.  by (3.3.10) in \cite{Elworthy-LeJan-Li-book},

 \begin{eqnarray*}
 {\hat D}U_t
 &=&\sum_{j=1}^{q+1}(-1)^{j+1}\rho(t)\langle T\xi_t(b^j), X(x_t)dB_t\rangle \cdot
 V_t^j\\
&&+\sum_{j=1}^{q+1}(-1)^{j+1}
\int_0^t \rho(s) \langle T\xi_s(b^j), X(x_s)dB_s\rangle \cdot
d\wedge^q\left(\breve \nabla X(-)dB_t)(V_t^j)\right)\\
&&+\sum_{j=1}^{q+1}(-1)^{j+1}
\int_0^t \rho(s) \langle T\xi_s(b^j), X(x_s)dB_s\rangle \cdot\\
&&\hskip 50pt
  \left(-{1\over 2}(\breve R^q)^*(V_t^j)dt
+d\wedge^q\breve \nabla A(V_t^j)dt\right)\\
&=&\sum_{j=1}^{q+1}(-1)^{j+1}\rho(t)\langle T\xi_t(b^j), X(x_t)dB_t\rangle \cdot
V_t^j\\
&&+ d\wedge^q\left(\breve \nabla X(-)dB_t)(U_t)\right)
-{1\over 2}(\breve R^q)^*(U_t)dt
+d\wedge^q\breve \nabla A(U_t)dt,
 \end{eqnarray*}
Taking the conditional expectation of the above equation,
as in the proof of Proposition 3.3.7 in \cite{Elworthy-LeJan-Li-book}
or that of Theorem A in \cite{Elworthy-Yor}, we have

\begin{eqnarray*}
 {\hat D}\overline{U_t}
&=&\overline{\sum_{j=1}^{q+1}(-1)^{j+1} \rho(t)
\langle T\xi_t(b^j), X(x_t)dB_t\rangle \cdot \wedge^q T\xi_t
 (b^1\wedge\dots \wedge\widehat{b^j}\wedge \dots \wedge b^{q+1})}\\
&&-{1\over 2}(\breve R^q)^*(\overline{U_t})dt
+d\wedge^q\breve \nabla A(\overline{U_t})dt\\
&=&\rho(t)\iota_{\langle-, X(x_t)dB_t\rangle}
   \overline{\wedge^{q+1} T\xi_t(b^1\wedge\dots \wedge b^{q+1})}
   -{1\over 2}(\breve R^q)^*(\overline{U_t})dt
+d\wedge^q\breve \nabla A(\overline{U_t})dt\\
&=&\rho(t)\iota_{\langle-,X(x_t)dB_t\rangle}
   \breve W_t^{A,q+1}(b^1\wedge\dots \wedge b^{q+1})
   -{1\over 2}(\breve R^q)^*(\overline{U_t})dt
+d\wedge^q\breve \nabla A(\overline{U_t})dt\\
&=&\rho(t)\iota_{\langle-,\breve{\paral_t} d\breve B_t\rangle}
   \breve W_t^{A,q+1}(b^1\wedge\dots \wedge b^{q+1})
   -{1\over 2}(\breve R^q)^*(\overline{U_t})dt
+d\wedge^q\breve \nabla A(\overline{U_t})dt
 \end{eqnarray*}
 Here $\breve{\paral_t}$ denotes parallel translation corresponding
 to the connection $\breve \nabla$ and $(\breve B_s)$ is the stochastic
 anti-development Brownian motion on $T_{x_0}M$, i.e. the martingale
 part of $\int_0^\cdot \breve{\paral_s}^{-1}\circ dx_s$.
Solve the equation to obtain
$$\overline{U_t}
  =\breve W_t^{A,q}\int_0^t \left( \breve W_s^{A,q}\right)^{-1}
   \rho(s) \iota_{\langle -,\breve {\parals_s} d\breve B_s\rangle}
   \breve W_s^{A,q+1}(b^1\wedge\dots \wedge b^{q+1}) \,ds.$$
Finally we arrive at:
\begin{corollary}
\label{cor-invariant}
Suppose that the torsion $\breve T\equiv XdY$ is invariant under the flow
 $\xi_t$. Then for $b^i\in T_{x_0}M, i=1,\dots, q+1$,
\begin{equation}
\begin{split}
&d(P_t\phi)(b^1\wedge\dots \wedge b^{q+1})\\
&=\left(\int_0^t\rho(r)dr\right)^{-1}
  \E\phi\left(\breve W_t^{A,q} \int_0^t \rho(s)(\breve W_s^{A,q})^{-1}
  \iota_{\langle -,\breve {\parals_s} d\breve B_s\rangle}
  \breve W_s^{A,q+1}(b^1\wedge\dots \wedge b^{q+1})\right)\\
 & -\phi \left( \breve W_t^{A,q}
     \left(\sum_{1\le i<j\le q+1}(-1)^{i+j+1} \breve T(b^i,b^j)\wedge
     b^1\wedge \dots \wedge\widehat{b^i} \dots
 \wedge\widehat{b^j} \dots \wedge b^{q+1}\right)\right).
\end{split}
\end{equation}
\end{corollary}

Note that in the non-invariant case, an analogous proof to that of Corollary
\ref{cor-invariant} leads to the intrinsic formula below.
 If $V$ is a $q+1$ vector, we define $\iota_{\breve T}V$
to be the operator from $\wedge^{q+1}TM \to \wedge^qTM$  which restricted
 to primitive vectors is given by:
$$\iota_{\breve T}(b^1\wedge \dots \wedge b^{q+1})=
\sum_{1\le i<j\le q+1}(-1)^{i+j+1}
     \breve T(b^i,b^j)\wedge
     b^1\wedge \dots \wedge\widehat{b^i} \dots
 \wedge\widehat{b^j} \dots \wedge b^{q+1}.$$

\begin{corollary}
Let $b$ be a $q+1$ vector in $\wedge^{q+1}T_{x_0}M$, then
\begin{equation}
\begin{split}
&\left(\int_0^t\rho(r)dr\right)\, d(P_t\phi)(b)\\
=&
  \E\phi\left(\breve W_t^{A,q} \int_0^t \rho(s)(\breve W_s^{A,q})^{-1}
 \left( \iota_{\langle -, \breve {\parals_s} d\breve B_s\rangle}
  \breve W_s^{A,q+1}(b)+\iota_{\breve T}   \breve W_s ^{A,q+1}(b)\,ds\right)\right).
\end{split}
\end{equation}
\end{corollary}

\noindent{\it Proof.}
We only need to worry about the last term of (\ref{Bismut-formula}), since
 the previous term is as in  Corollary \ref{cor-invariant}. For this
\begin{eqnarray*}
&&\sum_{1\le i< j\le q+1}(-1)^{i+j+1}\\
&&\E(\xi_t^*)(\phi)\left(\int_0^t\rho(s)T\xi_s^{-1}
 \left(\breve T\left(T\xi_s(b^i),T\xi_s(b^{j})\right)\right) ds
 \wedge b^1\wedge \dots \wedge\widehat{b^i} \dots
 \wedge\widehat{b^j} \dots \wedge b^{q+1}\right)\\
 &&=\E\phi\left(\sum_{1\le i<j\le q+1}(-1)^{i+j+1}
 T\xi_t\int_0^t\rho(s)T\xi_s^{-1}
 \left(\breve T\left(T\xi_s(b^i),T\xi_s(b^{j})\right)\right)ds \;
 \wedge \right.\\
&&\hskip 60pt \left. T\xi_t(b^1)\wedge \dots \wedge\widehat{T\xi_t(b^i)} \dots
 \wedge\widehat{T\xi_t(b^j)} \dots \wedge T\xi_t(b^{q+1})\right).
 \end{eqnarray*}
Set
\begin{eqnarray*}
U_t&=&\sum_{1\le i< j\le q+1}(-1)^{i+j+1}
 T\xi_t\int_0^t\rho(s)T\xi_s^{-1}
 \left(\breve T\left(T\xi_s(b^i),T\xi_s(b^{j})\right)\right)ds \;
 \wedge \\
&& T\xi_t(b^1)\wedge \dots \wedge\widehat{T\xi_t(b^i)} \dots
 \wedge\widehat{T\xi_t(b^j)} \dots \wedge T\xi_t(b^{q+1}).
 \end{eqnarray*}
Then after covariant differentiation and filtering we have
\begin{eqnarray*}
 {\hat D}\overline{U_t}
&=&\sum_{1\le i<j\le q+1}(-1)^{i+j+1}\rho(t)\cdot\\
&&\overline{ \breve T\left(T\xi_t(b^i),T\xi_t(b^j)\right)\wedge
 T\xi_t(b^1)\wedge \dots \wedge\widehat{T\xi_t(b^i)} \dots
 \wedge\widehat{T\xi_t(b^j)} \dots \wedge T\xi_t(b^{q+1})}\\
&&-{1\over 2}(\breve R^q)^*(\overline{U_t})dt
+d\wedge^q\breve \nabla A(\overline{U_t})dt\\
&=&\rho(t)\overline{\iota_{\breve T}
   \wedge^{q+1} T\xi_t(b^1\wedge\dots \wedge b^{q+1})}
   -{1\over 2}(\breve R^q)^*(\overline{U_t})dt
+d\wedge^q\breve \nabla A(\overline{U_t})dt\\
&=&\rho(t)\iota_{\breve T}
   \breve W_t^{A,q+1}(b)
   -{1\over 2}(\breve R^q)^*(\overline{U_t})dt
+d\wedge^q\breve \nabla A(\overline{U_t})dt,
 \end{eqnarray*}
 giving
 $$\bar U_t=\breve W_t^{A,q}\int_0^t \rho(s)\left(\breve W_s^{A,q}\right)^{-1}
   \iota_{\breve T}   \breve W_s ^{A,q+1}(b)\,ds.$$
The required equation now follows.   \hfill\rule{3mm}{3mm}

\subsubsection*{Special cases}
\begin{enumerate}
\item[(1)]
 When the connection $\breve \nabla$ defined by (\ref{7})  is the
Levi-Civita connection and $\rho(t)\equiv 1$, formula (\ref{Bismut-formula})
essentially reduces to (\ref{4}), but in this case $P_t$ has generator
 given by ${1\over 2}\Delta+\LL_A$ on smooth forms.

\item[(2)]
For a left invariant stochastic differential equation on a Lie group
$G$ with bi-invariant metric and $A\equiv 0$,
 formula  (\ref{Bismut-formula}) reduces to
\begin{equation}\label{formula-Li-group}
\begin{split}
&d(P_t\phi)(b^1\wedge\dots \wedge b^{q+1})\\
&=\left(\int_0^t\rho(r)dr\right)^{-1}
\E\left(\left\langle \int_0^t ad(x_s)dB_s, -\right\rangle \wedge
R_{x_t}^*(\phi)\right)(b^1\wedge\dots\wedge b^{q+1})\\
&-\sum_{1\le i<j\le q+1} (-1)^{i+j} \E\, R^*_{x_t}(\phi)
\left([b^i, b^j]\wedge b^1\wedge\dots \widehat{b^i}\wedge \dots
 \widehat{b^j}\wedge \dots b^{q+1}\right),
\end{split}
\end{equation}
for $x_0=e$ (so $b^j\in \g$, each $j$).
In this case the generator of $P_t$ is ${1\over 2}trace \nabla^R\nabla^R$
by (2.4.3) of \cite{Elworthy-LeJan-Li-book}.

\bigskip

 Formula  (\ref{formula-Li-group}) is intrinsic. It could have been deduced
from the path space integration by parts formula of \cite{Fang-Franchi97}.
The $q=0$ case was given in \cite{Elworthy-Li-JFA}.

\item[(3)] Another computable example
comes when $A\equiv 0$ and $X$ is chosen so that $\breve \nabla$ has
torsion $\breve T(u,v)={2\over n-1}(v\wedge u)Z(x)$ for
$u, v\in T_xM$ and $Z$ a fixed vector field on $M$. Here $v\wedge u$
is the operator such that
$(v\wedge u)Z(x)=\langle v,Z(x) \rangle_x u- \langle u, Z(x)\rangle_x v$.
This connection was used by \cite{Ikeda-Watanabe}. In
\cite{Elworthy-LeJan-Li-book}, example 2.3.5 (though there is a minor
misprint in the formula written there), the generator on $q$-forms
is shown to be
${1\over 2}\Delta+\LL_{[{2(q-1)\over n-1}-1]Z}-{2\over n-1}\iota_Zd$.
Set $Z^{\#}=\langle Z(x), -\rangle_x$. The term involving the torsion in
 (\ref{Bismut-formula}) reduces to
\begin{eqnarray*}
&&\left(\int_0^t\rho(r)dr\right)^{-1}{4\over n-1}
\E\,\xi_t^*(Z^{\#}\wedge \phi)(b^1\wedge\dots \wedge b^{q+1})\\
&&={4\over (n-1)\int_0^t\rho(r)\,dr}P_t(Z^{\#}\wedge\phi).
\end{eqnarray*}
In this case we have:
\begin{equation}
\label{41}
\begin{split}
&d(P_t\phi)(b^1\wedge \dots \wedge b^{q+1})\\
&=\left(\int_0^t\rho(r)dr\right)^{-1}
  \E\,\phi\left( \int_0^t \rho(s)  \iota_{\langle -,
  { \breve {\parals_s}} d \breve{B}_s\rangle}
  W_s^{q+1}(b)\right)\\
 & +  {4\over (n-1)\int_0^t\rho(r)\,dr}\left(Z^{\#}\wedge \phi\right)
  \left(b^1\wedge\dots\wedge b^{q+1}\right).
\end{split}
\end{equation}
\end{enumerate}

\bigskip

{\tiny \hskip-\parindent K. D. ELWORTHY, \hskip4pt MATHEMATICS
INSTITUTE,  WARWICK UNIVERSITY, COVENTRY CV4 7AL, UK}

{\tiny \hskip-\parindent XUE-MEI LI, \hskip4pt DEPARTMENT OF
COMPUTING AND MATHEMATICS, THE NOTTINGHAM TRENT UNIVERSITY, BURTON STREET,
NOTTINGHAM NG1 4BU, U.K.
and DEPARTMENT OF MATHEMATICS, UNIVERSITY OF CONNECTICUT,
196 AUDITORIUM ROAD, STORRS, CT 06269, USA.
} {\scriptsize e-mail:
xuemei.li@ntu.ac.uk}

\begin{thebibliography}{ELJL97}

\bibitem[Air76]{Airault76}
H{\'e}l{\`e}ne Airault.
\newblock Subordination de processus dans le fibr\'e tangent et formes
  harmoniques.
\newblock {\em C. R. Acad. Sci. Paris S\'er. A-B}, 282(22):Aiii, A1311--A1314,
  1976.

\bibitem[Bis81]{Bismut-Durham}
J.-M. Bismut.
\newblock Martingales, the {M}alliavin calculus and {H}\"ormander's theorem.
\newblock In {\em \em Stochastic integrals (Proc. Sympos., Durham, 1980),
  Lecture Notes in Mathematics 851}, pages 85--109. Springer, 1981.

\bibitem[Bis84]{Bismut-book}
J.~M. Bismut.
\newblock {\em Large deviations and the {M}alliavin calculus. {P}rogress in
  {M}ath. 45.}
\newblock Birkha\H user, 1984.

\bibitem[Dri92]{Driver92}
B.~K. Driver.
\newblock A {C}ameron-{M}artin type quasi-invariance theorem for {B}rownian
  motion on a compact {R}iemannian manifold.
\newblock {\em J. Functional Analysis}, 100:272--377, 1992.

\bibitem[DT01]{Driver-Thalmaier-2001}
B.~Driver and A.~Thalmaier.
\newblock Heat equation derivative formulas for vector bundles.
\newblock {\em J. Funct. Anal.}, 183:42--108, 2001.

\bibitem[EL]{Elworthy-Li-Hodge-2}
K.~D. Elworthy and Xue-Mei Li.
\newblock An ${L}^2$ theory for differential forms on path spaces.
\newblock In preparation.

\bibitem[EL94]{Elworthy-Li-JFA}
K.~D. Elworthy and Xue-Mei Li.
\newblock Formulae for the derivatives of heat semigroups.
\newblock {\em J. Funct. Anal.}, 125(1):252--286, 1994.

\bibitem[EL96]{Elworthy-Li-Ito}
K.~D. Elworthy and Xue-Mei Li.
\newblock A class of integration by parts formulae in stochastic analysis {I}.
\newblock In {\em \em It\^{o}'s Stochastic Calculus and Probability Theory
  (dedicated to It\^o on the occasion of his eightieth birthday)}. Springer,
  1996.

\bibitem[EL98]{Elworthy-Li-forms-CR}
K.~D. Elworthy and Xue-Mei Li.
\newblock {Bismut type formulae for differential forms}.
\newblock {\em C. R. Acad. Sci., S\'er/ I, Math. Paris}, 327(1):87--92, 1998.

\bibitem[EL00]{Elworthy-Li-Hodge-1}
K.~D. Elworthy and Xue-Mei Li.
\newblock Special {I}t{\^o} maps and an ${L}^2$ {H}odge theory for one forms on
  path spaces.
\newblock {\em Canadian Mathematical Society Con. Proc.}, 28:145--162, 2000.

\bibitem[ELJL97]{Elworthy-LeJan-Li-Tani}
K.~D. Elworthy, Y.~Le~Jan, and Xue-Mei Li.
\newblock Concerning the geometry of stochastic differential equations and
  stochastic flows.
\newblock In K.D. Elworthy, S.~Kusuoka, and I.~Shigekawa, editors, {\em New
  Trends in stochastic Analysis', Proc. Taniguchi Symposium, Sept. 1995,
  Charingworth}. World Scientific Press, 1997.

\bibitem[ELL99]{Elworthy-LeJan-Li-book}
K.~D. Elworthy, Y.~LeJan, and Xue-Mei Li.
\newblock {\em On the geometry of diffusion operators and stochastic flows,
  Lecture Notes in Mathematics 1720}.
\newblock Springer, 1999.

\bibitem[Elw82]{Elworthy-book}
K.~D. Elworthy.
\newblock {\em Stochastic {D}ifferential {E}quations on {M}anifolds, London
  Mathematical Society Lecture Notes Series 70}.
\newblock Cambridge University Press, 1982.

\bibitem[Elw88]{Elworthy-Stflour}
K.~D. Elworthy.
\newblock Geometric aspects of diffusions on manifolds.
\newblock In {\em P.~L. Hennequin, editor, Ecole d'Et\'e de Probabilit\'es de
  Saint-Flour XV-XVII, 1985-1987. Lecture Notes in Mathematics 1362, volume
  1362}, pages 276--425. Springer-Verlag, 1988.

\bibitem[Elw92]{Elworthy-flow}
K.~D. Elworthy.
\newblock Stochastic flows on {R}iemannian manifolds.
\newblock In M.~A. Pinsky and V.~Wihstutz, editors, {\em Diffusion processes
  and related problems in analysis, volume II. {B}irkhauser Progress in
  Probability}, pages 37--72. Birkhauser, Boston, 1992.

\bibitem[EY93]{Elworthy-Yor}
K.~D. Elworthy and M.~Yor.
\newblock Conditional expectations for derivatives of certain stochastic flows.
\newblock In J.~Az\'ema, P.A. Meyer, and M.~Yor, editors, {\em Sem. de Prob.
  XXVII. Lecture Notes in Mathematics 1557}, pages 159--172. Springer-Verlag,
  1993.

\bibitem[FF97]{Fang-Franchi97}
S.~Z. Fang and J.~Franchi.
\newblock De {R}ham-{H}odge-{K}odaira operator on loop groups.
\newblock {\em J. Functional Analysis}, 148:391--407, 1997.

\bibitem[IW81]{Ikeda-Watanabe}
N.~Ikeda and S.~Watanabe.
\newblock {\em Stochastic Differential Equations and Diffusion Processes}.
\newblock North-Holland, 1981.

\bibitem[JL91]{Jones-Leandre91}
J.~D.~S. Jones and R.~L{\'e}andre.
\newblock ${L}\sp p$-{C}hen forms on loop spaces.
\newblock In {\em Stochastic analysis (Durham, 1990), pages 103--162. London
  Mathematical Society Lecture Notes Series 167}. Cambridge University Press,
  Cambridge, 1991.

\bibitem[Kus88]{Kusuoka-degree}
S.~Kusuoka.
\newblock Degree theorem in certain {W}iener-{R}iemannian manifolds.
\newblock In {\em Stochastic {A}nalysis: {J}apanese-{F}rench {S}eminar 1987.
  Lecture Notes in Mathematics, 1322}, pages 93--108. Springer-Verlag, 1988.

\bibitem[Kus92]{Kusuoka-forms}
S.~Kusuoka.
\newblock Analysis on {W}iener spaces {II}, {D}ifferential forms.
\newblock {\em J. Funct. Anal.}, 103:229--274, 1992.

\bibitem[Lan62]{Lang-book}
S.~Lang.
\newblock {\em Introduction to differential manifolds}.
\newblock Interscience Publishers, 1962.

\bibitem[Lea]{leandre-survey}
R.~Leandre.
\newblock Analysis over loop space and topology.
\newblock To appear in Mathematical Notes.

\bibitem[Li92]{thesis}
Xue-Mei Li.
\newblock {\em Stochastic flows on noncompact manifolds}.
\newblock University of Warwick, 1992.
\newblock Ph.D. thesis.

\bibitem[Nor93]{Norris93}
J.~Norris.
\newblock Path integral formulae for heat kernels and their derivatives.
\newblock {\em Probability Theory and Related Fields}, 94:525--541, 1993.

\bibitem[RS84]{Rapoport-Sternberg84}
Diego Rapoport and Shlomo Sternberg.
\newblock On the interaction of spin and torsion.
\newblock {\em Ann. Physics}, 158(2), 1984.

\bibitem[Shi86]{Shigekawa-Hodge}
I.~Shigekawa.
\newblock De {R}ham-{H}odge-{K}odaira's decomposition on an abstract {W}iener
  space.
\newblock {\em J. Math. Kyoto Univ.}, 26(2):191--202, 1986.

\end{thebibliography}
\end{document}